\documentstyle[12pt]{article}
\begin{document}
\author{S. Ludkovsky, A. Khrennikov.}
\title{Stochastic processes on 
non-Archimedean spaces with values in non-Archimedean fields.}
\date{27 October 2001
\thanks{Mathematics subject classification
(1991 Revision) 28C20 and 46S10.} }
\maketitle
\begin{abstract}
Stochastic processes on topological vector spaces
over non-Archimedean fields and with
transition measures having values in non-Archimedean fields
are defined and investigated. For this the non-Archimedean analog 
of the Kolmogorov theorem is proved. The analogos of Markov
and Poisson processes are studied. For Poisson processes
the corresponding Poisson measures are considered
and the non-Archimedean analog of the L\`evy theorem is proved. 
Wide classes of stochastic processes are constructed. 
\end{abstract}
\section{Introduction.} 
Classical stochastic
analysis for real and complex
vector spaces and real or complex transition measures
is well developed 
\cite{dalfo,gihsko,gisk3,ikwat,mckean,malb,oksb},
but the stochastic analysis
on topological vector spaces over non-Archimedean fields and with
transition measures having values in non-Archimedean fields
was not studied. 
There are many differences of classical and non-Archimedean functional 
analysis \cite{boui,grus,hew,pietsch,put,reed,roo} and many theorems 
of classical functional analysis are not true in their classical form
in the non-Archimedean case, for example, measure theory,
operator theory, theory of function spaces.
This paper is devoted to such new non-Archimedean variant
of stochastic analysis and continues papers \cite{lustpr},
where real and complex valued transition measures 
of stochastic processes on non-Archimedean spaces were
considered. 
Stochastic differential equations on real Banach spaces and manifolds
are widely used for solutions of mathematical and physical
problems and for construction and investigation of measures
on them.
On the other hand, non-Archimedean functional analysis
develops fastly in recent years and also its applications
in mathematical physics 
\cite{ami,esc,sch1,roo,vla3,khrif,jang}.
Wide classes of quasi-invariant measures including 
analogous to Gaussian type
on non-Archimedean Banach spaces, loops and diffeomorphisms 
groups were investigated in \cite{lud,luumnls,luseamb,lubp99,lutmf99,lubp2}.
Quasi-invariant measures on topological groups and their
configuration spaces can be used 
for the investigations of their unitary representations
(see \cite{luseamb,lubp99,lutmf99,lubp2} and references therein).
\par In view of this developments non-Archimedean analogs
of stochastic equations and diffusion processes 
need to be investigated. Some steps in this direction were 
made in \cite{bikvol,evans}, where non-Archimedean time was
considered, but stochastic processes there were on spaces
of complex valued functions and transition measures were
real or complex valued. 
At the same time measures may be 
real, complex or with values in a non-Archimedean field.
\par In the classical stochastic analysis indefinite integrals
are widely used, but in the non-Archimedean case the field of $p$-adic
numbers $\bf Q_p$ has not linear order structure apart from $\bf R$. 
\par This work treats the case which was not considered by 
another authors. These investigations are not restricted
by the rigid geometry class \cite{freput}, since it is rather narrow.
Wider classes of functions and manifolds are considered. This is possible
with the use of
Schikhof's works on classes of functions $C^n$ in the sence of
difference quotients, which he investigated few years later the
published formalism of the rigid geometry.
\par Here are considered spaces of functions with values in Banach spaces
over non-Archimedean local fields, in particular, with values in
the field $\bf Q_p$ of $p$-adic numbers.
For this
non-Archimedean analogs of stochastic processes are considered
on spaces of functions with values in the non-Archimedean 
field such that a parameter 
analogous to the time is either real, $p$-adic 
or more generally can take values in any group (see \S \S 4.1, 4.2).
Certainly this encompasses cases of the time parameter
with values in adeles and ideles.
Their existence is proved in Theorem 4.3.
\par This became possible due to results of \S 2, where
the non-Archimedean variant of the Kolmogorov theorem 
was proved.
\par In \S 3 non-Archimedean analogs of Markov cylindrical distributions
are defined
and Propositions 3.3.1 and 3.3.2 about their boundedness
and unboundedness are proved. 
\par Poisson measures and processes play very important role
in classical stochastic analysis \cite{itkean,coxmil}.
In Section 5 their non-Archimedean analogs are considered.
All results of this paper are obtained for the first time.
\section{$p$-Adic probability measures.}  
Let $X$ be a set and $\cal R$ be a covering ring of $X$ such that
elements of $\cal R$ are subsets of $X$. Consider a field
$\bf K$ with a nontrivial non-Archimedean valuation such that 
${\bf K}\supset \bf Q_p$.Suppose that $\bf K$ is complete
as the ultrametric space.
\par {\bf 2.1. Definition.} Suppose that $\cal S$ is a subfamily of
$\cal R$ such that for each $A$ and $B$ in $\cal S$ there exists
$C\in \cal S$ with $C\subset A\cap B$, then $\cal S$ is called
shrinking. For a function $f: {\cal R} \to \bf K$ or
$f: {\cal R} \to \bf R$ the notation $\lim_{A\in \cal S}f(A)=0$
means that for each $\epsilon >0$ there exists $B\in \cal S$
such that $|f(A)|\le \epsilon $ for each $A\in \cal S$
with $A\subset B$. 
\par {\bf 2.2. Definition.} A mapping $\mu : {\cal R} \to \bf K$
is called a measure if it satisfies the following conditions:
\par $(i)$ $\mu (A\cup B)=\mu (A)+\mu (B)$
for each $A$ and $B$ in ${\cal R}$ such that $A\cap B=\emptyset $;
\par $(ii)$ for each $A\in \cal R$ its $\mu $-norm
$\| A\| _{\mu }:=\sup \{ |\mu (B)|: B\in {\cal R}, B\subset A \} <\infty $
is bounded;
\par $(iii)$ if ${\cal S}\subset \cal R$ is shrinking and
$\cap {\cal S}:=\bigcap_{S\in \cal S}=\emptyset $, then
$\lim_{A\in \cal S} \mu (A)=0$.
\par {\bf 2.3. Note.} These conditions are called respectively
additivity, boundedness and continuity. Condition $(iii)$ 
is equivalent to $\lim_{A\in \cal S} \| A\| _{\mu }=0$ for each
shrinking subfamily ${\cal S}$ in ${\cal R}$ with
$\cap {\cal S}=\emptyset $.
\par {\bf 2.4. Definition.} A measure $\mu : {\cal R}\to \bf K$
is called a probability measure if $\mu (X)=1$ and 
$\| X\| _{\mu }=:\| \mu \|=1$.
\par {\bf 2.5. Remarks.} For functions $f: X\to \bf K$
and $\phi : X\to [0,\infty )$ put $\| f\| _{\phi }:=\sup_{x\in X}
|f(x)|\phi (x)$. Consider the following function:
$$(1)\quad N_{\mu }(x):=\inf_{U: x\in U\in \cal R} \| U\| _{\mu }$$
for each $x\in X$. Put $\| f\| _{\mu }:=\| f\| _{N_{\mu }}$.
Then for each $A\subset X$ the function $\| A\|_{\mu }:=\sup_{
x\in A} N_{\mu }(x)$ is defined such that its restriction on
$\cal R$ coincides with that of given by Equation $2.2.(ii)$
(see also Chapter 7 \cite{roo}). 
A $\cal R$-step function $f$ is a function $f: X\to \bf K$
such that it is a finite linear combination over $\bf K$ of characteristic
functions $Ch_U$ of $U\in \cal R$. A function $f$ is called $\mu $-integrable
if there exists a sequence $\{ f_n: n\in {\bf N} \} $ of step functions
such that $\lim_{n\to \infty } \| f-f_n\| _{N_{\mu }} =0$. The Banach space
of $\mu $-integrable functions is denoted by $L(\mu ):=L(X,{\cal R},\mu ,{\bf K})$. 
There exists a ring ${\cal R}_{\mu }$ of subsets $A$ in $X$ for which
$Ch_A\in L(\mu )$. The ring ${\cal R}_{\mu }$ is the extension of the ring $\cal R$
such that ${\cal R}_{\mu }\supset \cal R$.
\par For example, if $\bf K$ is locally compact, then the valuation group
$\Gamma _{\bf K}:=\{ |x|: x\in {\bf K}, x\ne 0 \} $ is discrete in $(0,\infty )
\subset {\bf R}$. If $\mu $ is a measure such that
$0< \| \mu \| <\infty $, then there exists
$a\in \bf K$ such that $|a|=\| \mu \| ^{-1}$, since 
$\| \mu \| \in \Gamma _{\bf K}$ for discrete $\Gamma _{\bf K}$,
hence  $a\mu $ is also the measure with $\| \mu \|=1$.
If $\| \mu \|=1$, then $\mu $ is the nonzero measure.
For such $\mu $ with $\mu (X)=:b_X\in \bf K$
if $b_X\ne 1$ we can take new set $Y$ and define on $X_0:=
Y\cup X$ a minimal ring ${\cal R}_0$ generated by ${\cal R}$
and $\{ Y \} $, that is, ${\cal R}_0\cap Y= \{ \emptyset , \{ Y \} \} $
and ${\cal R}_0={\cal R}\cup \{ Y \} $. 
Since $\| \mu \|=1$, then $|b_X|\le 1$.
Put $\mu (Y):=1-b_X$, then there exists the extension
of $\mu $ from ${\cal R}$ on ${\cal R}_0$ such that
$\| \mu \|=1$ and $\mu (X_0)=1$, since $|1-b_X|\le \max (1,|b_X|)=1$. 
In particular, we can take a singleton $Y=\{ y \} $.
Therefore, probability measures are rather naturally related with
nonzero bounded measures.
This also shows that from $ \| \mu \| =1$ in general does not follow
$\mu (X)=1.$ Evidently, from $\mu (X)=1$ in general does not follow
$\| \mu \| =1$, for example, $X=\{ 0,1 \} $, ${\cal R}= \{ \emptyset ,
\{ 0 \}, \{ 1 \}, X \} $, $\mu (\{ 0\} )=a$, $\mu (\{ 1\} )=1-a$,
where $|a|>1$, hence $\| \mu \| =|a|>1.$
Wide class of probability $\bf Q_p$-valued measures on non-Archimedean 
Banach spaces was constructed in \S II.3.15 \cite{lulapm}.
\par Consider a nonvoid topological space $X$.
A topological space is called zero-dimensional if it has a base
of its topology consisting of clopen subsets.
A topological space $X$ is called a $T_0$-space
if for each two distinct points $x$ and $y$ in $X$
there exists an open subset $U$ in $X$ such that either
$x\in U$ and $y\in X\setminus U$ or
$y\in U$ and $x\in X\setminus U$.
\par A covering ring ${\cal R}$ of a space $X$ defines on
it a base of zero-dimensional topology $\tau _{\cal R}$
such that each element of ${\cal R}$
is considered as a clopen subset in $X$.
If $\pi : X\to Y$ is a mapping such that
$\pi ^{-1}({\cal R}_Y)\subset {\cal R}_X$, then a measure $\mu $
on $(X,{\cal R}_X)$ induces a measure $\nu :=\pi (\mu )$
on $(Y,{\cal R}_Y)$ such that $\nu (A)=\mu  (\pi ^{-1}(A))$
for each $A\in {\cal R}_Y$.
\par {\bf 2.6. Proposition.} {\it Let $(X,{\cal R},\mu )$ be a measure
space. Then there exists a quotient mapping $\pi : X\to Y$ on a Hausdorff
zero-dimensional space $(Y,\tau _{\cal G})$ and
$\pi (\mu ):=\nu $ is a measure on $Y$ such that
${\cal G}=\pi ({\cal R})$, where $(Y,{\cal G},\nu )$
is the measure space.}
\par {\bf Proof.} Suppose that $(Y,\tau _{\cal G})$
is a $T_0$-space, where $\cal G$ is a covering ring of $Y$.
For each two distinct points $x$ and $y$ in $Y$
there exists a clopen subset $U$ in $Y$ such that either
$x\in U$ and $y\in Y\setminus U$ or
$y\in U$ and $x\in Y\setminus U$, since the base of topology
$\tau _{\cal G}$ in $Y$ consists of clopen subsets.
On the other hand, $Y\setminus U$ is also clopen, since
$U$ is clopen. Therefore, $Y$ is the Hausdorff space.
Clearly this implies that $Y$ is the Tychonoff space (see \S 6.2
\cite{eng}, but it is necessary to note that we consider
the definition of the zero-dimensional space more general without
$T_1$-condition in \S 2.5).
\par Now we construct a $T_1$-space $Y$, that is a quotient space
of $X$. For this consider the relation in $X$: \\
$x\kappa y$ if and only if for each $S\in {\cal R}$ with $x\in S$
there is the inclusion $\{ x,y \} \subset S$.
Evidently, $x\kappa x$, that is, $\kappa $ is reflexive.
The relation $x\kappa y$ means, that $y\in V_x:=\bigcap_{x\in S\in {\cal R}}
S$, where $V_x$ is closed in $X$, then from $y\in S$ it follows, that
$x\in S$, since otherwise $y\notin V_x$, because ${\cal R}$
is a covering ring. Therefore, $V_x=V_y$
and $y\kappa x$, hence $\kappa $ is symmetric.
Let $x\kappa y$ and $y\kappa z$, then $V_x=V_y=V_z$, consequently,
$x\kappa z$ and $\kappa $ is transitive.
Therefore, $\kappa $ is the equivalence relation.
Let $\pi : X\to Y:=X/\kappa $ be the quotient mapping
and $Y$ be supplied with the zero-dimensional topology
generated by the covering ring ${\cal G}$ such that
$\pi ^{-1}({\cal G})={\cal R}$, since each $A\in \cal R$ is clopen,
then for each $x\in A\in \cal  R$ we have $V_x\subset A$.
Then $\pi ^{-1}([y])=V_y$ for each $y\in X$ and 
$[y]:=\pi (y)$. Hence each point $[y]\in Y$ is closed,
hence $Y$ is the $T_1$-space. The topology
in $Y$ is generated by the covering ring ${\cal G}$,
consequently, $Y$ is the Hausdorff space (see above), since
from the $T_1$ separation property it follows the $T_0$
separation property.
\par If ${\cal S}$ is the shrinking family with zero intersection
in $Y$ such that ${\cal S}\subset \cal G$, 
then $\pi ^{-1}({\cal S})$ is also the shrinking family
with zero intersection
in $X$ such that $\pi ^{-1}({\cal S})\subset \cal R$, 
hence from $\lim_{A\in \pi ^{-1}({\cal S})} \mu (A)=0$
it follows $\lim_{A\in {\cal S}} \nu (A)=0$.
Therefore, Condition $(iii)$ from \S 2.2 is satisfied.
Evidently, $\| \nu \| =\| \mu \|$
and $\nu $ is additive on ${\cal G}$, hence
$\nu $ is the measure.
\par {\bf 2.7.1. Note.} In view of Proposition 2.6
we consider henceforth Hausdorff zero-dimensional measurable
$(X,{\cal R})$ spaces if another is not specified.
\par In the classical case the principal role in stochastic analysis
plays the Kolmogorov theorem, that gives the possibility
to construct a stochastic process on the basis of a system of finite dimensional
(real-valued) probability distributions (see \S III.4 \cite{kolfp}).
The following three theorems resolve this problem for $\bf K$-valued
measures in cases of a product of measure spaces, a consistent family 
of measure spaces and in cases of bounded cylindrical distributions.
Finally Theorem 2.15.2 (the non-Archimedean analog of the Kolmogorov theorem)
as the particular case of Theorem 2.14 is formulated.
\par Consider now a family of probability measure
spaces $\{ (X_j, {\cal R}_j, \mu _j) : j \in \Lambda \} $,
where $\Lambda $ is a set.
Suppose that each covering ring ${\cal R}_j$ is complete relative
to the measure $\mu _j$, that is, ${\cal R}_j={\cal R}_{\mu _j}$, where
${\cal R}_{\mu _j}$ denotes the completion of ${\cal R}_j$
relative to $\mu _j$. Let $X:=\prod_{j\in \Lambda }
X_j$ be the product of topological spaces supplied with the
product (Tychonoff) topology $\tau _X$, where each $X_j$ is considered
in its $\tau _{{\cal R}_j}$-topology.
There is the natural continuous projection
$\pi _j: X\to X_j$ for each $j\in \Lambda $.
Let $\cal R$ be the ring of the form 
$\bigcup_{j_1,...,j_n\in \Lambda , n\in \bf N} 
\bigcap _{l=1}^n\pi _{j_l} ^{-1}({\cal R}_{j_l})$.
\par {\bf 2.7.2. Definition.} A triple
$(X,{\cal R},\mu )$ is called a cylindrical distribution
if it satisfies the following condition:
\par $\mu |_{\bigcap _{l=1}^n\pi _{j_l} ^{-1}({\cal R}_{j_l})}
=\prod_{l=1}^n\tilde \mu _{j_l}$ for each
$j_1,...,j_n\in \Lambda $ and $n\in \bf N$, where
$\tilde \mu _j(\pi _j^{-1}(A)):=\mu _j(A)$ for each
$A\in {\cal R}_j$; $\tilde \mu _j$ is the measure
on $(X,\pi _j^{-1}({\cal R}_j)).$
\par {\bf 2.8. Theorem.} {\it A cylindrical distribution
$\mu $ on $(X,{\cal R})$ has an extension
up to a probability measure $\mu $ on $(X,{\cal R}_{\mu })$, where
$\mu $ and $X$ are the same as in \S 2.7.}
\par {\bf Proof.} For each $j\in \Lambda $ we have the ring
${\cal R}_j$. Let $A$ and $B$ be in
$\bigcap _{l=1}^n\pi _{j_l} ^{-1}({\cal R}_{j_l})$, where
$j_1,..,j_n\in \Lambda $ and $n\in \bf N$.
Then $A=\bigcap _{l=1}^n\pi _{j_l} ^{-1}(A_l)$, where
$A_l\in {\cal R}_{j_l}$ for each $l=1,...,n$, analogously 
for $B$ with $B_l$ instead of $A_l$.
Such subsets $A$ form the base of the topology $\tau _X$
such that $\tau _X\supset \cal R$.
Therefore, $A\cap B$ and $A\setminus B$ and hence $A\cup B$
are in $\cal R$, since ${\cal R}_{j_1}\times ... \times
{\cal R}_{j_n}$ is the ring.
Therefore, $\cal R$ is the ring. 
The space $X$ in the topology $\tau _X$ is zero-dimensional,
since the base $\cal R$ of $\tau _X$ consists of clopen subsets
in $X$. It is necessary to verify that the triple $(X,{\cal R},\mu )$
satisfies Conditions $2.2.(i-iii).$ In general $\tau _X$ and $\cal R$
may not coincide, but as it is shown below the usage of the inclusion
$\tau _X\supset \cal R$ is sufficient for the proof.
On the other hand,
$(X_{j_1}\times ... \times X_{j_n}, {\cal R}_{j_1}\times ... \times
{\cal R}_{j_n}, \mu _{j_1}\times ... \times \mu _{j_n})$
is the measure space for each $j_1,...,j_n\in
\Lambda $ and $n\in \bf N$, consequently,
$\mu $ on $\cal R$ is additive.
For each $A$ of the outlined above form we have
\par $(i)$ $\| A \| _{\mu }=\prod_{l=1}^n\| A_l \| _{\mu _{j_l}}\le 1$.
Such elements $A$ in $\cal R$ form the base of Tychonoff topology
in $X$, consequently, $\| X\| _{\mu }\le 1$.
For each $j\in \Lambda $ we have $\mu _j(X_j)=1$, hence
$\mu (X)=1$ and $\| X\| _{\mu }=1$.
Therefore, $\mu $ satisfies Conditions $2.2.(i,ii)$.
For each $A\in \cal R$ the norm $\| A\| _{\mu }$ is defined.
\par Consider now the function $N_{\mu }(x)$ on
$(X,{\cal R})$, that is defined by the Formula 2.5.(1).
For each $\epsilon >0$ and $x\in X$ there exists
$A\in \cal R$ such that 
\par $(ii)$ $\| A\| _{\mu } -\epsilon <N_{\mu }(x)
\le \| A\| _{\mu }$. Each function $N_{\mu _j}(x_j)$
is upper semicontinuous on $(X_j,{\cal R}_j)$ by Theorem 7.6 \cite{roo}. In view of 
Lemma 7.2 \cite{roo} and Formula $(i)$
for each $x\in X$ and each $\epsilon >0$
there exists its neighborhood $A\in \cal R$ such that 
\par $(iii)$ $\prod _{l=1}^nN_{\mu _{j_l}}(y_{j_l})<
N_{\mu }(x)+\epsilon $
for each $y\in A$, where $y_j:=\pi _j(y)$ for each $j\in \Lambda $.
Hence for each $x\in X$ and each $\epsilon >0$
there exists its basic neighborhood $A$ such that
\par $(iv)$ $N_{\mu }(y)<N_{\mu }(x)+\epsilon $ for each $y\in A$, that is,
$N_{\mu }(x)$ is upper semicontinuous on $(X,{\cal R})$,
since $0\le N_{\mu _j}(x_j)\le 1$ for each $x_j\in X_j$ and $j\in \Lambda $.
From Formulas $(i,ii,iii)$ and $2.2.(ii)$ we have
\par $(v)$ $\| A\| _{\mu }=\sup_{x\in X}N_{\mu }(x)$ for each
$A\in \cal R$, \\
since $ \| A\|_{\mu }=\sup_{x_1\in A_1,...,x_n\in A_n}
\prod_{l=1}^nN_{\mu _{j_l}}(x_l)$.
\par Let $V$ be a compact subset of $X$. Then for each
$\epsilon >0$ it has a covering by clopen subsets
$E_x\in \cal R$ with $x\in V$ and $x\in E_x$
such that Inequalities
$(ii-iv)$ are satisfied for $E_x$ instead of $A$. In view of compactness of $V$
the covering $\{ E_x: x\in V \} $ has a finite subcovering
$\{ E_1,...,E_m \} $. Then $\bigcup_{l=1}^mE_l\in \cal R$.
Therefore,
\par $(vi)$ $\sup_{x\in V}N_{\mu }(x)\le
\max_{l=1,...,m} \| E_l\|_{\mu } \le \sup_{x\in V}N_{\mu }(x)+
2\epsilon $, \\
since $\sup_{x\in V}N_{\mu }(x)\le
\| \bigcup_{l=1}^mE_l \| _{\mu } $ due to Inequality $(ii)$, \\
$ \| \bigcup_{l=1}^mE_l \| _{\mu }=\max _{l=1,...,m}
\| E_l \| _{\mu }=\max _{l=1,...,m} \sup_{x\in E_l}N_{\mu }(x)$
due to Formula $(v)$ and Condition $2.2.(ii)$, since
$E_l\in \cal R$ for each $l=1,...,m$, but for each $E_l$ there 
exists $x_l\in E_l\cap V$ such that $N_{\mu }(y)
<N_{\mu }(x_l)+\epsilon $ for each $y\in E_l$ due to Inequality
$(iv)$ and the choise of $\{ E_1,...,E_m \} $ as the finite
subcovering of the covering $\{ E_x: x\in E_x\cap V, E_x\in {\cal R} \} $
(see above). Since $\epsilon >0$ is arbitrary
and $\| \bigcup_{l=1}^mE_l \|_{\mu } =\max_{l=1,...,m} \| E_l \|_{\mu } $, then 
\par $(vii)$ $\sup_{x\in V}N_{\mu }(x)=\inf_{{\cal R}\ni A\supset V}
\| A\|_{\mu }$, since for each $A\in \cal R$ 
such that $V\subset A$ there exists $\{ E_l: l=1,...,m \} $
with $m\in \bf N$, where each $E_l$ is as above, such that
$V\subset \bigcup_{l=1}^mE_l\subset A$. Though the compact subset $V$ 
is not necessarily in $\cal R$ we take
Equation $(vii)$ as the definition of $\| V\| _{\mu }:=
\inf_{{\cal R}\ni A\supset V}\| A\|_{\mu } $.
\par Now we verify, that $\mu $ on $(X,{\cal R})$ satisfies
Condition $2.2.(iii)$.
Let $\cal S$ be a shrinking subfamily in
$\cal R$ with $\bigcap {\cal S}=\emptyset $.
In view of Theorem 7.12 \cite{roo} for each $\epsilon >0$
the set $X_{j, \epsilon }:=\{ x: x\in X_j, N_{\mu _j}(x)\ge \epsilon \} $
is ${\cal R}_{\mu _j}$-compact.
For each $\delta >0$ choose a sequence $ \{ \epsilon _j:
\epsilon _j>0, j\in \Lambda \} $ with $\sup_{j\in\Lambda }
\epsilon _j<\delta $. In view of the Tychonoff theorem
(see \S 3.2.4 in \cite{eng})
$\prod_{j\in \Lambda }X_{j, \epsilon _j}=:X_{ \{ \epsilon _j:
j\} }$ is the compact subset in $X$. Since for each
$A$ and $B$ in $\cal S$ there exists $C\in \cal S$ such that
$C\subset A\cap B$ and $\cal R$ is the ring, then 
consider finite intersections of finite families in $\cal S$,
hence there exists
the minimal family ${\cal S}_0$ generated by $\cal S$
such that ${\cal S}_0\subset \cal R$,
${\cal S}\subset {\cal S}_0$ and ${\cal S}_0 $ is centered, that is,
$A\cap B\in {\cal S}_0 $ for each $A$ and $B$ in ${\cal S}_0 $.
Evidently, $\lim_{A\in \cal S} \| A\| _{\mu }=0$
is equivalent to $\lim_{A\in {\cal S}_0} \| A\| _{\mu }=0$, since
$\| B\| _{\mu }\le \| A\| _{\mu }$ for each $B\subset A$
(see also \S \S 2.3 and 2.5).
Denote ${\cal S}_0$ by $\cal S$ also.
\par Each element $S\in \cal S$ is clopen in the centered family
$\cal S$ and $X_{ \{ \epsilon _j: j\} }$ has the empty intersection
with $\bigcap \cal S$. In view of compactness of
$X_{ \{ \epsilon _j: j\} }$ there exists a finite subfamily
$S_1,...,S_n$ in $\cal S$ such that $X_{ \{ \epsilon _j: j\} }
\cap (\bigcap_{l=1}^nS_l)=\emptyset $, hence
$\lim_{A\in {\cal S}\cap X_{ \{ \epsilon _j: j\} }} \| A\| _{\mu }=0$, since
$\| \emptyset \| _{\mu }=0$. 
As above for each $\delta >0$ choose a finite covering $E_1,...,E_m
\in \cal R$ of $V:= X_{ \{ \epsilon _j: j\} }$ such that 
\par $(viii)$ $\| V\| _{\mu }\le \max_{l=1,...,m}\| E_l \| _{\mu }<
\| V\| _{\mu }+ \delta $ (see Equation $(vi)$), hence there exists
$A\in {\cal S} \cap ( \bigcup_{l=1}^mE_l ) $ such that $\| A\| _{\mu }
\le \delta $, where $( \bigcup_{l=1}^mE_l ) \in \cal R$.
\par For each $x\in X\setminus X_{ \{ \epsilon _j: j\} }$
there exists a basis neigbourhood $U=\bigcap_{l=1}^n
\pi ^{-1} (U_l)$ such that $U\cap X_{ \{ \epsilon _j: j\} }=
\emptyset $, since $X$ is Hausdorff and $X_{ \{ \epsilon _j: j\} }$
is compact, where $U_l\in {\cal R }_{j_l}$ and $j_l\in \Lambda $ for each $l=1,...,n$.
Therefore, 
\par $(ix)$ $N_{\mu }(x)\le \delta $ for each $x\in 
X\setminus X_{ \{ \epsilon _j: j\} }$, since $\| X_j \|_{\mu _j}=1$
and \\
$\sup_{x\in X\setminus X_{j,\epsilon _j}}N_{\mu _j}(x)
\le \epsilon _j<\delta $ for each $j\in \Lambda $. 
In view of Equations $(v,vi)$ we have
$\| A \| _{\mu } \le \delta $ for each $A\subset
X\setminus X_{ \{ \epsilon _j: j\} }$ such that $A\in \cal R$.
Then applying Equation $(viii)$ to $V=X_{ \{ \epsilon _j: j\} }$
we get 
$\| (\bigcap _{l=1}^nS_l)\cap (\bigcup_{k=1}^mE_k) \| _{\mu } \le \delta $
and $\| (\bigcap _{l=1}^nS_l)\cap (X\setminus (\bigcup_{k=1}^mE_k)) 
\| _{\mu } \le \delta $ due to Inequality $(ix)$,
where $X\setminus (\bigcup_{k=1}^mE_k)\in \cal R$, hence
$\| (\bigcap_{l=1}^nS_l)\| _{\mu } \le \delta $, since
$\| A\| _{\mu }\le \max ( \| A\cap B\|_{\mu },\| A\cap 
(X\setminus B)\|_{\mu })$, 
where $(\bigcap_{l=1}^nS_l)\in \cal S$. On the other hand,
$\delta >0$ is arbitrary, consequently,
$\lim_{A\in {\cal S}} \| A\| _{\mu }=0$. 
This means that $(X,{\cal R},\mu )$ is the measure space.
In view of Theorem 7.4 \cite {roo} the measure $\mu $ has an extension
$\mu $ on the completion ${\cal R}_{\mu }$ of $\cal R$ relative to
$\mu $, moreover, $\sup_{{\cal R}_{\mu }\ni B\subset A} |\mu (B)|=
\sup_{x\in A}N_{\mu }(x)$ for each $A\in {\cal R}_{\mu }$.
\par {\bf 2.9. Note.} Theorem 2.8 has an evident generalization
for bounded measures $\mu _j$ if two products $\prod_{j\in \Lambda _0}
\mu _j(X_j)\in \bf K$ and  $\prod_{j\in \Lambda } \| X_j \| _{\mu _j}   
<\infty $ converge, where $ \Lambda _0:= \{j:  j\in \Lambda ,
\mu _j(X_j)\ne 0 \} $, when $\Lambda \setminus \Lambda _0$ is finite. 
Since $\mu $ is defined on $\cal R$
and bounded on it, then $\mu $ has an extension to the bounded
measure $\mu $ on  ${\cal R}_{\mu }$ such that $\mu (X)=
\prod_{j\in \Lambda } \mu _j(X_j)$ and  $ \| X \| _{\mu } =
\prod_{j\in \Lambda } \| X_j \| _{\mu _j} $, where
${\cal R}_{\mu }$ is the completion of ${\cal R}$ relative to
$\mu $.
\par The conditions imposed in \S 2.7
on the family of measures are natural. 
In view of \S 2.8 if $\Lambda =\bf N$ and
$\prod_{j=1}^{\infty } \| X_j\| _{\mu _j}=0$, then
for each $x\in X$ and each $\epsilon >0$ there exists
$A\in \cal R$ such that $\prod_{j=1}^{\infty }\| \pi _j(A)\|_{\mu _j}
<\epsilon $, hence $\| X\| _{\mu }=0$. On the other hand, if 
$\prod_{j=1}^{\infty } \| X_j\| _{\mu _j}=\infty $, then
for each $r>0$ there exists $X _{ \{ \epsilon _j: j \} } $
and $A\supset X _{ \{ \epsilon _j: j \} }$ with
$A\in \cal R$ such that $\prod_{j=1}^{\infty }\| \pi _j(A)\|_{\mu _j}
>r$, hence $\| X\| _{\mu }=\infty $. 
If $\mu $ is a measure on $(X,{\cal R})$, then $\mu (A)\in \bf K$
for each $A\in \cal R$.
In the case of the product measure $\mu $ this leads to the restriction,
that $\prod_{j\in \Lambda _0}\mu _j(X_j)$ is convergent, when
$\Lambda \setminus \Lambda _0$ is finite. For infinite
$\Lambda \setminus \Lambda _0$ we have $\mu (A)=0$ for each
$A\in \cal R$.
The condition $\| X_j\| _{\mu _j}=1$ for each $j$ does not
guarantee this convergence. For example, if ${\bf K}=\bf Q_p$
with the prime $p>1$ and the set $\{ j: \mu _j(X_j)\in \{ 2,...,p-1 \} \} $
is infinite, then $\prod_{j\in \Lambda _0}\mu _j(X_j)$ diverges, since
the multiplicative group $({\bf Z}/{\bf pZ})\setminus \{ 0 \} $
of the quotient ring $({\bf Z}/{\bf pZ})$ is cyclic.
\par {\bf 2.10. Note.} A set $\Lambda $ is called directed
if there exists a relation $\le $ on it satisfying the following conditions:
\par $(D1)$ If $j\le k$ and $k\le m$, then $j\le m$;
\par $(D2)$ For every $j\in \Lambda $, $j\le j$;
\par $(D3)$ For each $j$ and $k$ in $\Lambda $ there exists
$m\in \Lambda $ such that $j\le m$ and $k\le m$.
A subset $\Upsilon $ of $\Lambda $ directed by $\le $ is called cofinal
in $\Lambda $ if for each $j\in \Lambda $ there exists
$m\in \Upsilon $ such that $j\le m$.
Suppose that $\Lambda $ is a directed set
and $\{ (X_j,{\cal R}_j,\mu _j ): j\in \Lambda \} $ 
is a family of probability measure spaces, where ${\cal R}_j$
is the covering ring (not necessarily separating). Supply each
$X_j$  with the topology $\tau _j$ such that its base is the ring
${\cal R}_j$ as in \S 2.5. 
Let this family be consistent in the following sence:
\par $(1)$ there exists a mapping $\pi ^k_j: X_k\to X_j$
for each $k\ge j$ in $\Lambda $ such that $(\pi ^k_j)^{-1}({\cal R}_j)
\subset {\cal R}_k$, $\pi ^j_j(x)=x$ for each $x\in X_j$ and each 
$j\in \Lambda $, $\pi ^m_k\circ \pi ^k_l=\pi ^m_l$ for each
$m\ge k\ge l$ in $\Lambda $;
\par $(2)$ $\pi ^k_l(\mu _k)=(\mu _l)$ for each $k\ge l$ in $\Lambda $.
Such family of measure spaces is called consistent.
\par {\bf 2.11. Theorem.} {\it Let 
$\{ (X_j,{\cal R}_j,\mu _j ): j\in \Lambda \} $ be a consistent family
as in \S 2.10. Then there exists a probability measure space
$(X,{\cal R}_{\mu },\mu )$ and a mapping $\pi _j: X\to X_j$ for each
$j\in \Lambda $ such that $(\pi _j)^{-1}({\cal R}_j)\subset {\cal R}$
and $\pi _j(\mu )=\mu _j$ for each $j\in \Lambda $.}
\par {\bf Proof.} We have $(\pi ^k_j)^{-1}({\cal R}_j)
\subset {\cal R}_k$ for each $k\ge j$ in $\Lambda $, then
$(\pi ^k_j)^{-1}(\tau _j)
\subset \tau _k$ for each $k\ge j$ in $\Lambda $, 
since each open subset in $(X_j,\tau _j)$
is the union of some subfamily $\cal G$ in ${\cal R}_j$
and $(\pi ^k_j)^{-1}(\bigcup {\cal G})=\bigcup_{A\in \cal G}
(\pi ^k_j)^{-1}(A)$. Therefore, each $\pi ^k_j$ is continuous
and there exists the inverse system
${\sf S}:=\{ X_k,\pi ^k_j,\Lambda \} $
of the spaces $X_k$. Its limit
$\lim {\sf S}=:X$ is the topological space with the topology
$\tau _X$ and continuous mappings $\pi _j: X\to X_j$
such that $\pi ^k_j\circ \pi _k= \pi _j$ for each $k\ge j$ in
$\Lambda $ (see \S 2.5 in \cite{eng}). Each element $x\in X$ is the thread
$x=\{ x_j: x_j\in X_j$ $\mbox{for each}$ $j\in \Lambda,$ 
$\pi ^k_j(x_k)=x_j$ $\mbox{for each}$ $k\ge j\in \Lambda \} $.
Then $\pi _j^{-1}({\cal R}_j)=:{\cal G}_j$ is the ring of subsets in $X$
such that ${\cal G}_j\subset \tau _X$ for each $j\in \Lambda $. 
The base of topology of $(X,\tau _X)$ is formed by subsets
$\pi _j^{-1}(A)$, where $A\in \tau _j$, $j\in \Lambda $, but ${\cal R}_j$
is the base of topology $\tau _j$ for each $j$, hence
$\{ B: B=\pi _j^{-1}(A),$ $A\in {\cal R}_j,$ $j\in \Lambda \} $
is the base of $\tau _X.$ Therefore, the ring ${\cal R}:=
\bigcup_{j\in \Lambda }{\cal G}_j$ is the base of $\tau _X.$
In view of Proposition 2.6 we can reduce our consideration
to the case when each ${\cal R}_j$ is separating on $X_j$
and ${\cal R}$ is separating on $X$, since
${\cal G}_j\subset {\cal G}_k$ for each $k\ge j$ in $\Lambda $.
\par Consider on $\cal R$ a function $\mu $ with values
in $\bf K$ such that $\mu (\pi _j^{-1}(A)):=\mu _j(A)$
for each $A\in {\cal R}_j$ and each $j\in \Lambda $.
If $A$ and $B$ are disjoint elements in $\cal R$, 
then there exists $j\in \Lambda $
such that $A$ and $B$ are in ${\cal G}_j$, hence 
\par $(i)$ $A=\pi _j^{-1}(C)$ and $B=\pi _j^{-1}(D)$ for some
$C$ and $D$ in ${\cal R}_j$, consequently,
$\mu (A\cup B)=\mu _j(C\cup D)=\mu _j(C)+\mu _j(D)=\mu _j(A)+\mu _j(B)$,
that is, $\mu $ is additive. Moreover, $\| A\|_{\mu }=\| C\|_{\mu _j}$
for each $A=\pi _j^{-1}(C)$ with $C\in {\cal R}_j$, hence
$\| X\| _{\mu }=1$. Since $\mu (X)=\mu _j(X_j)$ and $\mu _j(X_j)=1$
for each $j\in \Lambda $, then  $\mu (X)=1$. Therefore,
$\mu $ satisfies Conditions $2.2.(i,ii)$.
It remains to verify Condition $2.2.(iii)$.
By Formula $2.5.(1)$ we have the function $N_{\mu }(x)$
on $(X,{\cal R})$ such that for each $x\in X$ and $\epsilon >0$
there exists $A\in \cal R$ such that
\par $(ii)$ $\| A\|_{\mu }-\epsilon <N_{\mu }(x)\le  \| A\|_{\mu }$.
In view of $(i)$ and upper semicontinuity of $N_{\mu _j}(x_j)$
on $(X_j,{\cal R}_j)$ for each $x\in X$ and $\epsilon >0$
there exists $j\in \Lambda $ and its neighborhood $A=\pi _j^{-1}(C)\in \cal R$
such that 
\par $(iii)$ $N_{\mu _j}(y_j)<N_{\mu }(x)+\epsilon $
for each $y\in A$, where $y_j:=\pi _j(y)$.
Hence for each $x\in X$ and each $\epsilon >0$
there exists its basic neighborhood $A$ such that
\par $(iv)$ $N_{\mu }(y)<N_{\mu }(x)+\epsilon $ for each $y\in A$, that is,
$N_{\mu }(x)$ is upper semicontinuous on $(X,{\cal R})$,
since $0\le N_{\mu _j}(x_j)\le 1$ for each $x_j\in X_j$ and $j\in \Lambda $.
From Formulas $(i,ii,iii)$ and $2.2.(ii)$ we have
\par $(v)$ $\| A\| _{\mu }=\sup_{x\in X}N_{\mu }(x)$ for each
$A\in \cal R$, since $ \| A\|_{\mu }=\sup_{x\in C}
N_{\mu _j}(x)$. For a compact subset $V$ in $X$ and each $\epsilon >0$
there exists a finite covering $\{ E_1,...,E_m\} \subset \cal R$
of $V$ such that inequalities $(ii-iv)$ are satisfied for each 
$E_l$ instead of $A$.
Therefore,
\par $(vi)$ $\sup_{x\in V}N_{\mu }(x)\le
\max_{l=1,...,m} \| E_l\|_{\mu } \le \sup_{x\in V}N_{\mu }(x)+
2\epsilon $ and 
\par $(vii)$ $\sup_{x\in V}N_{\mu }(x)=\inf_{{\cal R}\ni A\supset V}
\| A\|_{\mu }$. Though the compact subset $V$ 
is not necessarily in $\cal R$ we take
Equation $(vii)$ as the definition of $\| V\| _{\mu }:=
\inf_{{\cal R}\ni A\supset V}\| A\|_{\mu } $.
\par Choose a sequence
$\epsilon _j=\delta >0$ for each $j\in \Lambda $, where
$\delta >0$ is independent from $j$.
For each $\epsilon _j>0$ a subset $X_{j,\epsilon _j}:=
\{ x_j: x_j\in X_j, N_{\mu _j}(x_j)\ge \epsilon _j>0 \} $ is compact.
If $x_k\in X_{k,\epsilon _k}$, then $N_{\mu _j}(\pi ^k_j(x_k))\ge
\epsilon _k$ for each $j<k$, since $(\pi ^k_j)^{-1}({\cal R}_j)\subset 
{\cal R}_k$ and $\| B\| _{\mu _k}\le \| A\| _{\mu _k}$
for each $B$ and $A$ in ${\cal R}_k$ with $B\subset A$.
Hence $\pi ^k_j( X_{k,\epsilon _k})
\subset X_{j,\epsilon _k}$ for each $j\le k$
in $\Lambda $. Since
$\pi ^m_k\circ \pi ^k_l=\pi ^m_l$ for each $m\ge k\ge l$
in $\Lambda $, then $\{ X_{k,\delta }, \pi ^k_j,
\Lambda \} $ is the inverse system. 
The image $\pi ^k_j(X_{k,\delta })$ of each compact set
$X_{k,\delta }$ is compact for each $k>j$
(see Theorem 3.1.10 
\cite{eng}), since each $(X_k,\tau _k)$ is the Hausdorff space
in our consideration.
Since the limit of an inverse mapping system
of compact spaces is compact (see Theorem 3.2.13 \cite{eng}), 
then the limit
$X_{ \{ \epsilon _j:j \} }:=\lim \{ X_{k,\epsilon _k}, \pi ^k_j,
\Lambda \} $ is the compact subset in $X$
such that $X_{ \{ \epsilon _j:j \} }$ is homeomorphic
with $\theta (X)\cap \prod_{k\in \Lambda }  X_{k,\epsilon _k},$
where $\theta : X\hookrightarrow \prod_{k\in \Lambda }  X_k$
is the embedding.
For a shrinking family $\cal S$ in $\cal R$ consider all
finite intersections of finite families in $\cal S$, this gives
a centered family  ${\cal S}_0$ in $\cal R$ 
such that  ${\cal S}\subset {\cal S}_0$
and denote  it also by $\cal S$.
Applying $(i-vii)$ to $V=X_{ \{ \epsilon _j:j \} }$ and using
basic neighborhoods $U=\pi _k^{-1}(U_k)$, where $U_k\in {\cal R}_k$,
we get analogously to \S 2.8 
that for each shrinking family $\cal S$ in $\cal R$
with $\bigcap {\cal S}=\emptyset $ there exists
$\lim _{A\in \cal S} \| A \| _{\mu }=0$, since
due to $(vi)$ we have 
\par $(viii)$ $N_{\mu }(x)\le \delta $ for each $x\in X\setminus 
X_{ \{ \epsilon _j:j \} }$, since
$N_{\mu _j}(x_j)< \delta $ for each $x_j\in X\setminus
X_{j,\delta }$ and each $j\in \Lambda $.
Using the completion of $\cal R$ relative to
$\mu $ we get the probability measure space
$(X,{\cal R}_{\mu },\mu )$.
\par {\bf 2.12. Remark.} Each family
$\{ {\cal G}_j: j\in \Lambda \}  $ such that for each
$j$ and $k$ in $\Lambda $ there exists $m$ such that
${\cal G}_m\supset {\cal G}_k\cup {\cal G}_j$ defines on the set
$\Lambda $ the structure of the directed set:
$k\ge j$ if and only if ${\cal G}_k\supset {\cal G}_j$.
If there exists a (continuous) retraction
$r$ of $\prod _{j\in \Lambda }
X_j=:Y$ on $\lim \{ X_k,\pi ^k_j,\Lambda \}=:X$, that is,
$r(Y)=X$ and $r(x)=x$ for each $x\in X$, then a measure
$\nu $ on $(Y,{\cal R}_Y)$ induces a measure $\mu =r(\nu )$
on $(X,{\cal R}_X)$, since $r^{-1}({\cal R}_X)\subset
{\cal R}_Y$, such that in this particular case
Theorem 2.11 follows from Theorem 2.8.
On the other hand, Theorem 2.8 can be deduced from Theorem 2.11,
since a product of topological spaces is the particular case
of a limit of an inverse system, but the direct proof
of \S 2.8 is simpler.
\par {\bf 2.13. Note.} Theorem 2.11 has an evident generalization
to the following case:
$ \| X_j \| _{\mu _j} <\infty $ for each $j$ and 
there exist two limits $\lim_{j\in \Lambda _0}
\mu _j(X_j)\in \bf K$ and  $\lim_{j\in \Lambda } \| X_j \| _{\mu _j}   
<\infty $, where $ \Lambda _0:= \{j:  j\in \Lambda ,
\mu _j(X_j)\ne 0 \} $ and $\Lambda \setminus \Lambda _0$
is bounded in $\Lambda $.
We have $ \| X_j \| _{\mu _j}\le \| X_k \| _{\mu _k}$
for each $j\le k$ in $\Lambda $, since $\pi ^k_j(\mu _k)=\mu _j$
and $(\pi ^k_j)^{-1}({\cal R}_k)\subset {\cal R}_j$.
Since $\Lambda $ is directed, then  
$\lim_{j\in \Lambda } \| X_j \| _{\mu _j}=\sup_{
j\in \Lambda } \| X_j \| _{\mu _j}$.
Since a cylindrical distribution
$\mu $ is defined on $\cal R$
and bounded on it, then $\mu $ has an extension to the bounded
measure $\mu $ on  ${\cal R}_{\mu }$ such that $\mu (X)=
\lim_{j\in \Lambda } \mu _j(X_j)$ and  $ \| X \| _{\mu } =
\lim_{j\in \Lambda } \| X_j \| _{\mu _j} $, where
${\cal R}_{\mu }$ is the completion of ${\cal R}$ relative to
$\mu $.
\par Let now $X$ be a set with a covering ring 
$\cal R$ such that $X\in \cal R$. 
Let also $\{ (X,{\cal G}_j,\mu _j): j\in \Lambda \} $
be a family of measure spaces such that  $\Lambda $
is directed and  ${\cal G}_j\subset {\cal G}_k$
for each $j\le k\in \Lambda $, ${\cal R}=\bigcup_{
j\in  \Lambda } {\cal G}_j$. Suppose $\mu : {\cal R}\to
\bf K$ is such that $\mu |_{{\cal G}_j}=\mu _j$ and
$\mu _k|_{{\cal G}_j}=\mu _j$ for each $j\le k$ in $\Lambda $.
Then the triple $(X,{\cal R},\mu )$ is called the cylindrical distribution.
For each $A\in \cal R$ there exists $j\in \Lambda $
such that $A\in {\cal G}_j$, hence $\| A\|_{\mu _j}=
\| A\|_{\mu _k}$ for each $k\ge j$ in $\Lambda $, consequently,
$\| A\|_{\mu }:=\lim_{k\in \Lambda }\| A\| _{\mu _k}$ is correctly defined.
Suppose $\mu $ is bounded, that is,$\| X\|_{\mu }<\infty $.
(A particular simpler case is given below in \S 2.15).
\par {\bf 2.14. Theorem.} {\it Let $(X,{\cal R},\mu )$ 
be a bounded cylindrical distribution as in \S 2.13. Then $\mu $ 
has an extension to a bounded 
measure $\mu $ on the completion ${\cal R}_{\mu }$
of $\cal R$ relative to $\mu $.}
\par {\bf Proof.} Let $\tau _X$ be a topology on $X$
generated by the base $\cal R$. In view of Proposition 2.6
each covering ring ${\cal G}_j$ of $X$ produces
an equivalence relation $\kappa _j$ and a quotient mapping
$\pi _j: X\to X_j$ such that $\pi _j({\cal G}_j)=:{\cal R}_j$
is a separating covering ring of $X_j$, where $X_j$ is zero-dimensional
and Hausdorff. Moreover, ${\cal R}_j$ is the base of topology
$\tau _j$ on $X_j$. Since ${\cal G}_k\supset {\cal G}_j$
for each $k\ge j$, then on $(X_k,(\pi ^k_j)^{-1}({\cal R}_j))$
there exists an equivalence relation $\kappa ^k_j$ and a quotient
(continuous) mapping $\pi ^k_j: X_k\to X_j$ such that
$\pi ^m_k\circ \pi ^k_j=\pi ^m_j$ for each $j\le k\le m$ in
$\Lambda $. Hence there exists an inverse mapping system
$\{ X_k,\pi ^k_j,\Lambda \} $.
Therefore, the set $X$ in the topology $\tau _X$
generated by its base $\cal R$ consisting of clopen subsets
is homeomorphic with $\lim \{ X_k,\pi ^k_j,\Lambda \} $.
Each $\pi _j(\mu )=\mu _j$ is a bounded measure on $(X_j,{\cal R}_j)$
such that $\pi ^k_j(\mu _k)=\mu _j$ and $(\pi ^k_j)^{-1}({\cal R}_j)\subset
{\cal R}_k$ for each $k\ge j\in \Lambda $. Therefore,
$\{ (X_j,{\cal R}_j,\mu _j): j\in \Lambda \} $
is the consistent family of measure spaces.
From the definition of $\mu $ it follows that $\mu $ is additive,
hence $\| X\|_{\mu }$ is correctly defined. From 
$\| X\|_{\mu }<\infty $ it follows $\| X_j\|_{\mu _j}<\infty $
for each $j\in \Lambda $  and there exists $\lim_{j\in \Lambda }
\| X_j\|_{\mu _j}=\| X\|_{\mu }$. From $X\in \cal R$ it follows, that
$\mu (X)=\mu _j(X_j)$ for each $j\in \Lambda $.
Then this Theorem follows from Theorem 2.11.
\par {\bf 2.15.1. Note.} Let $X:=\prod_{t\in T}X_t$
be a product of sets $X_t$ and on $X$ a covering ring
$\cal R$ be given such that for each $n\in \bf N$ and pairwise
distinct points $t_1,...,t_n$ in a set $T$ there exists a
$\bf K$-valued measure $P_{t_1,..,t_n}$ on a covering ring
${\cal R}_{t_1,...,t_n}$ of $X_{t_1}\times ... \times X_{t_n}$
such that $\pi ^{t_1,...,t_{n+1}}_{t_1,...,t_n}({\cal R}_{t_1,...,
t_{n+1}})={\cal R}_{t_1,...,t_n}$ for each $t_{n+1}\in T$
and $P_{t_1,...,t_{n+1}}(A_1\times ... \times A_n\times X_{n+1})=
P_{t_1,...,t_n}(A_1\times ... \times A_n)$
for each $A_1\times ... \times A_n\in {\cal R}_{t_1,...,t_n}$,
where $\pi ^{t_1,...,t_{n+1}}_{t_1,...,t_n}: X_{t_1}\times ...
\times X_{t_{n+1}}\to  X_{t_1}\times ... \times X_{t_n}$
is the natural projection, $A_l\subset X_{t_l}$ for each
$l=1,...,n$. Suppose that the cylindrical distribution
is bounded, that is, 
\par $\sup_{t_1,...,t_n\in T, n\in \bf N}
\| P_{t_1,...,t_n}\| <\infty $ and there exists 
\par $\lim_{t_1,...,t_n \in T_0; n\in \bf N} 
P_{t_1,...,t_n}(X_{t_1}\times ... \times X_{t_n})\in \bf K$,
where $T_0:=\{ t\in T: P_t(X_t)\ne 0 \} $, $T\setminus T_0$ is finite.
\par {\bf 2.15.2. Theorem} (the non-Archimedean analog of the 
Kolmogorov theorem){\bf .} {\it The cylindrical distribution 
$P_{t_1,...,t_n}$ from \S 2.15.1 has an extension
to a bounded measure $P$ on the completion ${\cal R}_P$
of ${\cal R}:=\bigcup_{t_1,...,t_n\in T, n\in \bf N}
{\cal G}_{t_1,...,t_n}$ relative to $P$, where ${\cal G}_{t_1,...,t_n}:=
(\pi _{t_1,...,t_n})^{-1}({\cal R}_{t_1,...,t_n})$
and $\pi _{t_1,...,t_n}: X\to X_{t_1}\times ... \times X_{t_n}$
is the natural projection.}
\section{Markov distributions for a non-Archimedean Banach space.}
\par {\bf 3.1. Remark.} Let $H=c_0(\alpha ,{\bf K})$ be a Banach space over
a non-Archimedean field $\bf K$ with an ordinal $\alpha $(that is useful due to
Kuratowski-Zorn lemma, see \cite{eng,roo})
and the standard orthonormal base $\{ e_j: j\in \alpha \} $,
$e_j=(0,...,0,1,0,...)$ with $1$ on the $j$-th place, that is,
$c_0(\alpha ,{\bf K})=\{ x:$ $x=(x_i: i\in \alpha , x_i\in {\bf K} )$,
$\sup_{i\in \alpha }|x_i|=:\| x\| <\infty ,$ for each $b>0$ a set
$\{ i: |x_i|>b \} $ is finite $ \} $.
Suppose that $\bf K$ is complete as the ultrametric space.
For example, $\bf K$ is such that ${\bf Q_p}\subset {\bf K}\subset
\bf C_p$ or ${\bf F_p(\theta )}\subset {\bf K}$, where
$p$ is a prime number, $\bf Q_p$ is the field
of $p$-adic numbers, $\bf C_p$ is the field of complex
$p$-adic numbers, ${\bf F_p(\theta )}$ is the field
of formal power series by an indeterminate $\bf \theta $ 
over the finite field $\bf F_p$ consisting of $p$ elements.
Let $\sf U^P$ be a cylindrical ring generated by projections 
$\pi _F: H\to F$ on
finite dimensional over $\bf K$ subspaces $F$ in $H$ and  
rings $Bco(F)$ of clopen subsets. 
This ring ${\sf U^P}$ is the base of the weak topology $\tau _{H,w}$
in $H$. 
Each vector $x\in H$ is considered as continuous linear functional 
on $H$ by the formula $x(y)=\sum_jx^jy^j$ for each $y\in H$, 
so there is the natural embedding $H\hookrightarrow H^*=l^{\infty }
(\alpha ,{\bf K})$, where $x=\sum_jx^je_j$, $x^j\in \bf K$,
$l^{\infty }(\alpha ,{\bf K}):=
\{ x:$ $x=(x_i: i\in \alpha , x_i\in {\bf K} )$,
$\sup_{i\in \alpha }|x_i|=:\| x\| <\infty \} $.
This justifies the following generalization.
\par {\bf 3.2. Notes and definitions.} Let $T$ 
be a subset in $\bf \Lambda $ and containing a point $t_0$
and $X_t=X$ be a locally $\bf K$-convex space
for each $t\in T$, where $\bf \Lambda $ is an additive group, for example,
$\bf \Lambda $ is contained in $\bf R$ or $\bf C$ or a 
non-Archimedean field. Put $(\tilde X_T,
{\tilde {\sf U}}):=\prod_{t\in T}(X_t,{\sf U}_t)$ be a product of
measurable spaces, where ${\sf U}_t$ are rings 
of clopen subsets of $X_t$, $\tilde {\sf U}$ is the ring of cylindrical
subsets of $\tilde X_T$ generated by projections $\tilde \pi _q: \tilde X_t\to
X^q$, $X^q:=\prod_{t\in q}X_t,$ $q\subset T$ is a finite subset of $T$
(see \S I.1.3 \cite{dalfo}).
Let $\bf K_s$ be a subfield of $\bf C_s$such that
$\bf K_s$ is complete as the ultrametric space, where $s$
is a prime number.
A function $P(t_1,x_1,t_2,A)$ with values in $\bf K_s$ 
for each $t_1\ne t_2\in T$, $x_1\in X_{t_1}$, $A\in {\sf U}_{t_2}$
is called a transition measure if it satisfies the following conditions:\\
$(i)\mbox{ the set function }\nu _{x_1,t_1,t_2}(A)
:=P(t_1,x_1,t_2,A)\mbox{ is 
a } \mbox{measure on }$ $(X_{t_2},{\sf U}_{t_2});$\\
$(ii)\mbox{ the function }\alpha _{t_1,t_2,A}(x_1)
:=P(t_1,x_1,t_2,A)\mbox{ of the variable } x_1$ $\mbox{ is  }
{\sf U}_{t_1}-\mbox{measurable},$ that is,
$\alpha _{t_1,t_2,A}^{-1}(Bco({\bf K_s}))\subset {\sf U}_{t_1}$; \\
$$(iii)\mbox{ }P(t_1,x_1,t_2,A)=\int_{X_z}
P(t_1,x_1,z,dy)P(z,y,t_2,A)\mbox{ for each }t_1\ne t_2\in T,$$
that is, $P(z,y,t_2,A)$ as the function by $y$
is in $L((X_z,{\sf U}_z),\nu _{x_1,t_1,z},{\bf K_s})$.
A transition measure  $P(t_1,x_1,t_2,A)$ is called normalised if
$$(iv)\mbox{ }P(t_1,x_1,t_2,X_{t_2})=1\mbox{ for each }t_1\ne t_2\in T.$$
For each set $q=(t_0,t_1,..,t_{n+1})$ of pairwise distinct points in $T$
there is defined a measure in $X^g:=\prod_{t\in g}X_t$ by the formula
$$(v)\mbox{ }\mu ^q_{x_0}(E)=\int_E\prod_{k=1}^{n+1}P(t_{k-1},x_{k-1},
t_k,dx_k), \mbox{ }E\in {\sf U}^g:=\prod_{t\in g}{\sf U}_t,$$
where $g=q\setminus \{ t_0 \}$, variables $x_1,...,x_{n+1}$ are such that
$(x_1,...,x_{n+1})\in E$, $x_0\in X_{t_0}$ is fixed.
\par Let $E=E_1\times X_{t_j}\times E_2$, where 
$E_1\in \prod_{i=1}^{j-1}{\sf U}_{t_i}$, $E_2\in \prod_{i=j+1}^{n+1}
{\sf U}_{t_i}$, if the transition measure 
$P(t,x_1,t_2,dx_2)$ is normalised, then
\par $(vi)$ $\mu ^q_{x_0}(E)=\int_{E_1\times E_2}[\prod_{k=1}^{j-1}P(
t_{k-1},x_{k-1},t_k,dx_k)]\times $ 
$[\int_{X_{t_j}}P(t_{j-1},x_{j-1},t_j,dx_j)$
$\prod_{k=j+1}^{n+1}P(t_{k-1},x_{k-1},t_k,dx_k)]$ 
$=\mu ^r_{x_0}(E_1\times E_2),$ where $r=q\setminus \{ t_j \} .$
From Equation $(vi)$ it follows, that 
$$(vii)\mbox{ }[\mu ^q_{x_0}]^{\pi ^q_v}=\mu ^v_{x_0}$$ 
for each $v<q$ (that is, $v\subset q$), where $\pi ^q_v: X^g\to X^w$
is the natural projection, $g=q\setminus \{ t_0 \} ,$
$w=v\setminus \{ t_0 \} .$ 
Therefore, due to Conditions $(iv,v,vii):$ $\{ \mu ^q_{x_0}; \pi ^q_v;
\Upsilon _T \} $ is the consistent family of measures, which induce
the cylindrical distribution
$\tilde \mu _{x_0}$ on $(\tilde X_T, \tilde {\sf U})$
such that $\tilde \mu _{x_0}(\pi
_q^{-1}(E))=\mu ^q_{x_0}(E)$ for each $E\in {\sf U}^g$, where 
$\Upsilon _T$ is the family of all finite subsets $q$ in $T$
such that $t_0\in q\subset T$, $v\le q\in \Upsilon _T$, $\pi _q: \tilde X_T\to
X^g$ is the natural projection, $g=q\setminus \{ t_0 \} $.
\par The cylindrical distributions given by Equations $(i-v,vii)$
are called Markov distributions (with time $t\in T$).
\par {\bf 3.3.} {\bf Proposition. 1.} {\it If a normalized transition measure
$P$ satisfies the condition
$$(i)\mbox{ }C:=\sup_q[\sum_{k=1}^nln (\sup_x
\| \nu _{x,t_{k-1},t_k} \|)] < \infty ,$$ where $q=(t_0,t_1,...,t_n)$
with pairwise distinct points $t_0,..,t_n \in T$ and $n\in \bf N$,
then the Markov cylindrical distribution ${\tilde \mu }_{x_0}$ is bounded
and it has an extension to a bounded measure ${\tilde \mu }_{x_0}$
on the completion $\tilde {\sf U}_{{\tilde \mu }_{x_0}}$ of
$\tilde {\sf U}$ relative to ${\tilde \mu }_{x_0}$.}
\par {\bf 3.3.2. Proposition.} {\it If 
$$(ii)\mbox{ }C_x:=\sup_q[\sum_{k=1}^nln 
\| \nu _{x,t_{k-1},t_k} \|]=\infty $$ 
for each $x,$ where $q=(t_0,t_1,...,t_n)$
with pairwise distinct points $t_0,..,t_n \in T$ and $n\in \bf N$,
then the Markov cylindrical distribution
${\tilde \mu }_{x_0}$ has the unbounded variation
on each nonvoid set $E\in \tilde {\sf U}.$}
\par {\bf Proof.} (1). If $E\in \tilde {\sf U}$, 
then $E\in {\sf U}^g$ for some set
$q=(t_0,t_1,...,t_n)$ with pairwise distinct points
$t_0,...,t_n \in T$ and $n\in \bf N$ and $g=q\setminus \{ t_0 \} $,
consequently,$|\mu ^q_{x_0}(E)|\le \prod_{k=1}^n
\sup_x\| \nu _{x,t_{k-1},t_k}\|$ 
$\le exp(C)<\infty ,$ since $t_k\in T$ for each $k=0,1,...,n$, hence
$\sup_{q,E} |\mu ^q_{x_0}(E)|=\| \tilde \mu _{x_0} \| \le exp(C)$.
In view of Theorem 2.15.2 we get an extension of $\tilde \mu _{x_0}$
to a bounded measure on $\tilde {\sf U}_{{\tilde \mu }_{x_0}}$.
\par (2). For each $(t_1,t_2,x)$ with $x$ in $\pi _{t_0,t_2}(E)$
there exists a set
$\delta (t_1,t_2,x)\in {\sf U}_{t_2}\cap \pi _{t_0,t_2}(E)$ such that
$\| \delta (t_1,t_2,x)\| _{\nu _{x_1,t_1,t_2}}>1+\epsilon (t_1,t_2,x_1,x)$,
where $\epsilon (t_1,t_2,x_1,x)>0$. In view of Condition $(ii)$ 
for each $R>0$ and $x$ we choose
$q$ such that $\sum_{k=1}^n\epsilon (t_k,t_{k+1},x_1,x)>R$. 
For chosen $u\ne u_1\in T$ and $x\in \pi _{t_0,u}(E)\subset
X_u$ we represent the set $\delta (u,u_1,x)$
as a finite union of disjoint subsets $\gamma _{j_1}
\in {\sf U}_{u_1}$ such that for each
$\gamma _{j_1}$ and $u_2\ne u_1$ there is a set $\delta _{j_1}
\in {\sf U}_{u_2}\cap \pi _{t_0,u_2}(E)$
satisfying $\| \delta _{j_1}\| _{\nu _{x_1,u_1,u_2}}
\ge 1+\epsilon (u_1,u_2,x_1,x)$for each $x\in \gamma _{j_1}.$
Then by induction $\delta _{j_1,...,j_n}=\bigcup_{j_{n+1}=1}^{m_{n+1}}
\gamma _{j_1,...,j_{n+1}}$ so that for $u_{n+2}\ne u_{n+1}\in T$
there is a set $\delta _{j_1,...,j_{n+1}}
\in {\sf U}_{u_{n+1}}\cap \pi _{t_0,u_{n+1}}(E)$ for which
$\| \delta _{j_1,...,j_{n+1}} \| _{\nu _{x_{n+1},u_{n+1},u_{n+2}}}
\ge 1+\epsilon (u_{n+1},u_{n+2},x_{n+1},x)$ for each $x\in 
\gamma _{j_1,...,j_{n+1}}.$
Put $\Gamma ^{u,x_0}_{j_1,...,j_n}=$ $\{ x: x(u)=x_0, x(u_1)
\in \gamma _{j_1},...,x(u_n)\in \delta _{j_1,...,j_n},$
$x(u_{n+1})\in \gamma _{j_1,...,j_n} \} $
and $\Gamma ^{u,x_0}:=(\bigcup_{j_1,...,j_n}
{\Gamma ^{u,x_0}}_{j_1,...,j_n})\in \tilde {\sf U}$, since
$m_1\in \bf N$,...,$m_n\in \bf N$. 
Then \\
$\| \Gamma ^{u,x_0}\| _{\tilde \mu _{x_0}}\ge
\sup_{j_1,...,j_n}\| (\prod_{k=1}^{n+1} \nu _{u_{k-1},
x_{k-1},u_k}(dx_k)|_{ \gamma _{j_1}\times ...\times \gamma _{j_1,...,j_n}
\times \delta _{j_1,...,j_n}} \| $ \\
$\ge \prod_{k=1}^n[1+\epsilon (u_{k-1},u_k,x_{k-1},x_k)]>R$,
consequently, $\| E\| _{\tilde \mu _{x_0}}=\infty ,$ \\
since $\| E\| _{\tilde \mu _{x_0}}\ge \sup_{\Gamma ^{u,x_0}}
\| \Gamma ^{u,x_0}\| _{\tilde \mu _{x_0}}$ and $R>0$ is arbitrary.
\par {\bf 3.4.} Let $X_t=X$ for each $t\in T$, ${\tilde X}_{t_0,x_0}:=
\{ x\in {\tilde X}_T: $ $x(t_0)=x_0 \} .$ We define a projection
operator $\bar \pi _q:$ $x\mapsto x_q$, where $x_q$ is defined on 
$q=(t_0,...,t_{n+1})$ such that $x_q(t)=x(t)$ for each $t\in q$,
that is, $x_q=x|_q$. For every $F: {\tilde X}_T\to \bf C_s$
there corresponds $(S_qF)(x):=F(x_q)=F_q(y_0,...,y_n),$ where $y_j=x(t_j)$,
$F_q: X^q\to \bf C_s$. We put ${\sf F}:=$ $\{ F| F: {\tilde X}_T\to {\bf C_s},$
$S_qF\mbox{ are }{\sf U}^q-\mbox{measurable} \} $.
If $F\in \sf F$, $\tau =t_0\in q$, then there exists an integral
$$(i)\mbox{ }J_q(F)=\int_{X^q}(S_qF)(x_0,...,x_n)\prod_{k=1}^{n+1}P(t_{k-1},
x_{k-1},t_k,dx_k).$$
\par {\bf Definition.} A function $F$ is called integrable with respect 
to the Markov cylindrical distribution $\mu _{x_0}$ if the limit
$$(ii)\mbox{ }\lim_qJ_q(F)=:J(F)$$ along the generalized net by 
finite subsets $q$ of $T$ exists. This limit is called a functional 
integral with respect to the Markov cylindrical distribution:
$$(iii)\mbox{ }J(F)=\int_{{\tilde X}_{t_0,x_0}}F(x)\mu _{x_0}(dx).$$
\par {\bf 3.5. Remark.} Consider a $\bf K_s$-valued measure 
$P(t,A)$ on $(X,{\sf U})$ for each $t\in T$
such that $A-x\in \sf U$ for each $A\in \sf U$ and $x\in X$, 
where $A\in \sf U$, $X$ is a locally $\bf K$-convex space,
$\sf U$ is a covering ring of $X$.
Suppose $P$ be a spatially homogeneous transition measure (see also \S 3.2), 
that is,
$$(i)\mbox{ }P(t_1,x_1,t_2,A)=P(t_2-t_1,A-x_1)$$ 
for each $A\in \sf U$, $t_1\ne t_2 \in T$ 
and $t_2-t_1\in T$ and every $x_1\in X$, where $P(t,A)$ 
satisfies the following condition:
$$(ii)\mbox{ }P(t_1+t_2,A)=\int_XP(t_1,dy)P(t_2,A-y)$$
for each $t_1$ and $t_2$ and $t_1+t_2\in T$.
Such a transition measure $P(t_1,x_1,t_2,A)$ is 
called homogeneous. In particular for $T=\bf Z_p$
we have 
$$(iii)\mbox{ }P(t+1,A)=\int_XP(t,dy)P(1,A-y).$$ 
If $P(t,A)$ is a continuous function by $t\in T$ 
for each fixed $A\in \sf U$, then Equation $(iii)$ defines
$P(t,A)$ for each $t\in T$, when $P(1,A)$ is given, since
$\bf Z$ is dense in $\bf Z_p$.
\S \S 2.7 and 2.15 and 4.3 provide examples of Markov distributions.
Examples of Markov distributions are also Poisson and Gaussian
distributions given below and in a forthcoming paper.
\par {\bf 3.6. Notes and definition.}
 Let $X$ be a locally $\bf K$-convex space
and $P$ satisfies Conditions $3.2(i-iii)$. 
For $x$ and $z \in \bf Q_p^n$ we denote by
$(z,x)$ the following sum $\sum_{j=1}^n x_jz_j$, where $x=(x_j:$
$j=1,...,n)$, $x_j \in \bf Q_p$. Each number $y\in \bf Q_p$ has  a
decomposition $y=\sum_l
a_lp^l$, where $a_l \in (0,1,...,p-1)$,
$\min (l:$ $a_l\ne 0)=:ord_p(y)> - \infty $  for
$y\ne 0$ and $ord(0):=\infty $ \cite{nari,sch1}.
We define a symbol $\{ y \}_p:=\sum_{l<0} a_lp^l$ for
$|y|_p>1$ and $\{ y\} _p=0$ for $|y|_p \le 1$. 
We consider a character 
of $X$, $\chi _{\gamma }: X\to {\bf C_s}$  given by the following formula:
$$(i)\mbox{ }\chi _{\gamma }(x)=\epsilon ^{z^{-1}\{ (e,\gamma (x)) \} _p}$$ 
for each $ \{ (e,\gamma (x)) \} _p\ne 0$,
$\chi _{\gamma }(x):=1$ for $ \{ (e,\gamma (x)) \} _p=0,$
where $\epsilon =1^z$ is a root of unity, $z=p^{ord(\{ (e,\gamma (x)) \} _p)},$ 
$\gamma \in X^*$, $X^*$ denotes the topologically 
conjugated space of continuous $\bf K$-linear functionals
on $X$, the field $\bf K$ as the $\bf Q_p$-linear space is
$n$-dimensional, that is,
$dim_{\bf Q_p}{\bf K}=n$, $\bf K$ as the Banach space over $\bf Q_p$
is isomorphic with $\bf Q_p^n$, $e=(1,...,1)\in \bf Q_p^n$,
where $s\ne p$ are prime numbers (see \cite{vla3} and \cite{lulapm}). Then 
$$(ii)\mbox{ }\phi (t_1,x_1,t_2,y):=\int_X\chi _y(x)P(t_1,x_1,t_2,dx)$$
is the characteristic functional of the transition measure $P(t_1,x_1,t_2,dx)$
for each $t_1\ne t_2\in T$ and each $x_1\in X$.
In the particular case of $P$ satisfying Conditions $3.5.(i,ii)$ with $t_0=0$
its characteristic functional is such that 
$$(iii)\quad \phi (t_1,x_1,t_2,y)=\psi (t_2-t_1,y)\chi _y(x_1),\mbox{ where}$$
$$(iv)\quad \psi (t,y):=\int_X\chi _y(x)P(t,dx)\mbox{ and}$$
$$(v)\quad \psi (t_1+t_2,y)=\psi (t_1,y)\psi (t_2,y)$$ 
for each $t_1\ne t_2\in T$ and $t_2-t_1\in T$ and
$t_1+t_2\in T$ respectively and $y\in X^*$, $x_1\in X$.
\section{Non-Archimedean stochastic processes.}
\par {\bf 4.1. Remark and definition.} A measurable space $(\Omega ,{\sf F})$ 
with a probability $\bf K_s$-valued
measure $\lambda $ on a covering ring $\sf F$ of a set
$\Omega $  is called a probability space and is denoted by
$(\Omega ,{\sf F},\lambda )$. Points $\omega \in \Omega $ are called 
elementary events and values $\lambda (S)$  
probabilities of events $S\in \sf F$. A measurable map
$\xi : (\Omega ,{\sf F})\to (X,{\sf B})$ is called a random variable
with values in $X$, where 
${\sf B}$ is a covering ring such that
${\sf B}\subset Bco(X)$, $Bco(X)$ is the ring
of all clopen subsets of a locally $\bf K$-convex space $X$,
$\xi ^{-1}({\sf B})\subset \sf F$,
where $\bf K$ is a non-Archimedean field complete as an ultrametric space. 
\par The random variable $\xi $ induces a normalized measure $\nu _{\xi }(A):=
\lambda (\xi ^{-1}(A))$ in $X$ and a new probability space 
$(X,{\sf B},\nu _{\xi }).$
\par Let $T$ be a set with a covering ring $\cal R$ and a measure
$\eta : {\cal R}\to \bf K_s$. Consider the following Banach space
$L^q(T,{\cal R},\eta ,H)$ as the completion of the set
of all ${\cal R}$-step functions $f: T\to H$ relative to the following 
norm:
\par $(1)\quad \| f\|_{\eta ,q}:=\sup_{t\in T}\| f(t)\|_H
N_{\eta }(t)^{1/q}$ for $1\le q<\infty $ and
\par $(2)\quad \| f\|_{\eta ,\infty }:=\sup_{1\le q<\infty }
\| f(t)\|_{\eta ,q}$, where $H$ is a Banach space over $\bf K$
(see also \S 2.5). For $0<q<1$ this is the metric space with the metric
\par $(3)\quad \rho _q(f,g):=\sup_{t\in T}\| f(t)-g(t)\|_H
N_{\eta }(t)^{1/q}.$  
\par If $H$ is a complete locally $\bf K$-convex space,
then $H$ is a projective limit of Banach spaces
$H=\lim \{ H_{\alpha },\pi ^{\alpha }_{\beta }, \Upsilon \} $,
where $\Upsilon  $ is a directed set, $\pi ^{\alpha }_{\beta }:
H_{\alpha }\to H_{\beta }$ is a $\bf K$-linear continuous mapping
for each $\alpha \ge \beta $, $\pi _{\alpha }: H\to H_{\alpha }$
is a $\bf K$-linear continuous mapping such that $\pi ^{\alpha }_{\beta }
\circ \pi _{\alpha }=\pi _{\beta }$ for each $\alpha \ge \beta $
(see \S 6.205 \cite{nari}).
Each norm $p_{\alpha }$ on $H_{\alpha }$ induces a prednorm
${\tilde p}_{\alpha }$ on $H$. If $f: T\to H$, then $\pi _{\alpha }\circ
f=:f_{\alpha }: T\to H_{\alpha }$. In this case $L^q(T,{\cal R},\eta ,H)$
is defined as a completion of a family of all step functions
$f: T\to H$ relative to the family of prednorms
\par $(1')\quad \| f\|_{\eta ,q,\alpha }:=\sup_{t\in T}{\tilde p}_{\alpha }
(f(t))N_{\eta }(t)^{1/q}$, $\alpha \in \Upsilon $, for $1\le q<\infty $ and
\par $(2')\quad \| f\|_{\eta ,\infty ,\alpha  }:=\sup_{1\le q<\infty }
\| f(t)\|_{\eta ,q,\alpha }$, $\alpha \in \Upsilon $,
or pseudometrics
\par $(3')\quad \rho _{q,\alpha }(f,g):=\sup_{t\in T}{\tilde p}_{\alpha }
(f(t)-g(t))N_{\eta }(t)^{1/q}$, $\alpha \in \Upsilon $, for $0<q<1$.
Therefore, $L^q(T,{\cal R},\eta ,H)$ is isomorphic with the projective limit \\
$\lim \{ L^q(T,{\cal R},\eta ,H_{\alpha }),\pi ^{\alpha }_{\beta },
\Upsilon \} $.
For $q=1$ we write simply $L(T,{\cal R},\eta ,H)$ and
$\| f\|_{\eta }$. This definition is correct, since
$\lim_{q\to \infty }a^{1/q}=1$ for each $\infty >a>0$.
For example, $T$ may be a subset of $\bf R$. Let $\bf R_d$
be the field $\bf R$ supplied with the discrete topology. Since
the cardinality $card ({\bf R})={\sf c}=2^{\aleph _0},$
then there are bijective mappings of $\bf R$ on $Y_1:=\{ 0,...,b \}^{\bf N}$
and also on $Y_2:={\bf N}^{\bf N}$, where $b$ is a positive integer number.
Supply $\{ 0,...,b \}$ and $\bf N$ with the discrete topologies
and $Y_1$ and $Y_2$ with the product topologies.
Then zero-dimensional spaces
$Y_1$ and $Y_2$ supply $\bf R$ with covering separating rings
${\cal R}_1$ and ${\cal R}_2$ contained in $Bco(Y_1)$ and $Bco(Y_2)$
respectively. Certainly this is not related
with the standard (Euclidean) metric in $\bf R$.
Therefore, for the space $L^q(T,{\cal R},\eta ,H)$ we can consider
$t\in T$ as the real time parameter. If $T\subset \bf F$
with a non-Archimedean field $\bf F$, then we can consider
the non-Archimedean time parameter.
If $T$ is a zero-dimensional $T_1$-space, then denote by
$C^0_b(T,H)$ the Banach space of continuous bounded functions
$f: T\to H$ supplied with the norm:
\par $(4)\quad \| f\|_{C^0}:=\sup_{t\in T} \| f(t)\|_H<\infty $. \\
If $T$ is compact, then $C^0_b(T,H)$ is isomorphic with
the space $C^0(T,H)$ of continuous functions. 
\par For a set $T$ and a complete locally $\bf K$-convex 
space $H$ over $\bf K$ consider the 
product $\bf K$-convex space $H^T:=\prod_{t\in T}H_t$ in the product topology,
where $H_t:=H$ for each $t\in T$.
\par Then take on either $X:=X(T,H)=L^q(T,{\cal R},\eta ,H)$ or $X:=
X(T,H)=C^0_b(T,H)$ or on $X=X(T,H)=H^T$ a covering ring ${\sf B}$ such that
${\sf B}\subset Bco(X)$. Consider a random variable
$\xi : \omega \mapsto \xi (t,\omega )$ with values in $(X,{\sf B})$
and $t\in T$.
\par Events $S_1,...,S_n$ are called independent in total if
$P(\prod_{k=1}^nS_k)=\prod_{k=1}^nP(S_k)$. Subrings
${\sf F}_k\subset {\sf F}$ are said to be independent if
all collections of events $S_k\in {\sf F}_k$ are independent in total, 
where $k=1,...,n$, $n\in \bf N$. To each collection of random variables
$\xi _{\gamma }$ on $(\Omega ,{\sf F})$ with $\gamma \in \Upsilon $
is related the minimal ring ${\sf F}_{\Upsilon }\subset \sf F$
with respect to which all $\xi _{\gamma }$ are measurable, where $\Upsilon $
is a set.
Collections $\{ \xi _{\gamma }: $ $\gamma \in \Upsilon _j \} $
are called independent if such are 
${\sf F}_{\Upsilon _j}$, where $\Upsilon _j\subset \Upsilon $ for each
$j=1,...,n,$ $n\in \bf N$.
\par Consider $T$ such that $card(T)>n$. For $X=C^0_b(T,H)$ or $X=H^T$
define $X(T,H;(t_1,...,t_n);(z_1,...,z_n))$ as a closed submanifold
of $f: T\to H$, $f\in X$ such that $f(t_1)=z_1,...,f(t_n)=z_n$, where
$t_1,...,t_n$ are pairwise distinct points in $T$ and
$z_1,...,z_n$ are points in $H$.
For pairwise distinct points $t_1,...,t_n$ in $T$ with
$N_{\eta }(t_1)>0,...,N_{\eta }(t_n)>0$ define 
$X(T,H;(t_1,...,t_n);(z_1,...,z_n))$ as a closed submanifold
which is the completion relative to the norm $\| f\|_{\eta ,q}$
of a family of $\cal R$-step functions $f: T\to H$ such that
$f(t_1)=z_1,...,f(t_n)=z_n$. In these cases 
$X(T,H;(t_1,...,t_n);(0,...,0))$ is the proper $\bf K$-linear subspace 
of $X(T,H)$ such that $X(T,H)$ is isomorphic with
$X(T,H;(t_1,...,t_n);(0,...,0))\oplus H^n$, since if
$f\in X$, then $f(t)-f(t_1)=:g(t)\in X(T,H;t_1;z_1)$
(in the third case we use that $T\in \cal R$ and hence there exists
the embedding $H\hookrightarrow X$). For $n=1$ and $t_0\in T$
and $z_1=0$ we denote $X_0:=X_0(T,H):=X(T,H;t_0;0)$.
\par {\bf 4.2. Defintion.} We define a (non-Archimedean)
stochastic process $w(t,\omega )$ with values in $H$ 
as a random variable such that:
\par $(i)$ the differences $w(t_4,\omega )-w(t_3,\omega )$ 
and $w(t_2,\omega )-w(t_1,\omega )$ are independent
for each chosen $\omega $, $(t_1,t_2)$ and $(t_3,t_4)$ with $t_1\ne t_2$,
$t_3\ne t_4$, either $t_1$ or $t_2$ is not in the two-element set
$ \{ t_3,t_4 \} ,$ where $\omega \in \Omega ;$
\par $(ii)$ the random variable $\omega (t,\omega )-\omega (u,\omega )$ has 
a distribution $\mu ^{F_{t,u}},$ where $\mu $ is a probability 
$\bf K_s$-valued measure on $(X(T,H),{\sf B})$ from \S 4.1, 
$\mu ^g(A):=\mu (g^{-1}(A))$ for $g: X\to H$ such that
$g^{-1}({\cal R}_H)\subset \sf B$ and each $A\in
{\cal R}_H$, a continuous linear operator $F_{t,u}: X\to H$ is given by
the formula $F_{t,u}(w):=w(t,\omega )-w(u,\omega )$ 
for each $w\in L^q(\Omega ,{\sf F},\lambda ;X_0),$
where $1\le q\le \infty ,$ $X_0$ 
is the closed subspace of $X$ as in \S 4.1,
${\cal R}_H$ is a covering ring of $H$ such that
$F_{t,u}^{-1}({\cal R}_H)\subset \sf B$ for each $t\ne u$ in $T$;
\par $(iii)$ we also put $w(0,\omega )=0,$  
that is, we consider a $\bf K$-linear subspace
$L^q(\Omega ,{\sf F},\lambda ;X_0)$
of $L^q(\Omega ,{\sf F},\lambda ;X)$,
where $\Omega \ne \emptyset $.
\par It is seen that $w(t,\omega )$ is a Markov process with 
transition measure $P(u,x,t,A)=\mu ^{F_{t,u}}(A-x)$.
\par This definiton is justified by the following Theorem.
\par {\bf 4.3.} {\bf Theorem.} {\it 
Let either $X=C^0_b(T,H)$ or $X=H^T$ or $X=L^g(T,{\cal R},\eta ,H)$ with
$1\le g\le \infty $ be the same spaces as in \S 4.1,
where the valuation group $\Gamma _{\bf K}$ is discrete
in $(0,\infty )$.
Then there exists a family $\Psi $
of pairwise inequivalent (non-Archimedean) stochastic processes on 
$X$of the cardinality $card (\Psi )\ge card (T) card (H)$
or $card (\Psi )\ge card ({\cal R}) card (H)$
respectively.}
\par {\bf Proof.} Each complete locally $\bf K$-convex space $H$ is 
a projective limit of Banach spaces $H_{\alpha }$. Therefore,
due to \S \S 2.5 and 4.2 it is sufficient to consider the case of
the Banach space $H$. Since $H$ is over the field 
$\bf K$ with the valuation group $\Gamma _{\bf K}$ discrete
in $(0,\infty )$, then $H$ is isomorphic with the Banach space
$c_0(\alpha ,{\bf K})$ (see Theorems 5.13 and 5.16
\cite{roo}), where $\alpha $ is an ordinal.
\par Let ${\cal R}_{\bf K}$ be a covering separating ring
of $\bf K$ such that elements of ${\cal R}_{\bf K}$ are
clopen subsets in $\bf K$. Then there exists a lot of probability measures
$m$ on $({\bf K},{\cal R}_{\bf K})$ with values in $\bf K_s$,
for example, atomic measure with atoms $a_j$ such that for each
$U\in {\cal R}_{\bf K}$ either $a_j\subset U$ or $a_j\subset
{\bf K}\setminus U$. For example, this can be done for singleton atoms.
Let the family $\Upsilon $ of $a_j$ be countable and $\lim_j m(a_j)=0$, when
$\Upsilon $ is infinite. Then $m(S):=\sum_{a_j\subset S}a_j$
for each $S\in {\cal R}_{\bf K}$ and $\| m\| =\sup_j|m(a_j)|$.
If $\bf K$ is infinite and contains a locally compact infinite
subfield $\bf F$ with a nontrivial valuation, then $\bf K$ can be
considered as a locally $\bf F$-convex space.
As the locally $\bf F$-convex space $\bf K$ in its weak topology 
is isomorphic with ${\bf F}^{\gamma }$, since $\Gamma _{\bf F}$
is discrete in $(0,\infty )$ and there is the non-Archimedean variant
of the Hahn-Banach theorem, where $\gamma $ is a set (see \cite{nari,roo}).
Having a measure on $\bf F$ we can construct a probability measure
on $\bf K$ due to Theorem 2.8 and Note 2.9 and Remarks 2.5.
\par Therefore, consider also the particular case of the locally compact field
$\bf K$. If $\bf K$ is infinite, then either ${\bf K}\supset \bf Q_p$
or ${\bf K}=\bf F_p(\theta )$ with the corresponding prime number
$p$, since $\bf K$ is with the nontrivial valuation \cite{wei}.
If $\bf K$ is finite, then ${\bf K}=\bf F_p$. 
Let $s$ be a prime number such that $s\ne p$, then $\bf K$
is $s$-free as the additive topological group (see the Monna-Springer theorem
in \S 8.4 \cite{roo}).
Therefore, there exists the $\bf K_s$-valued Haar measure
$w$ on $\bf K$, that is, the bounded measure on each clopen compact subset
of $\bf K$ with $w(B({\bf K},0,1))=1$
and $w(y+A)=w(A)$ for each $A\in Bco({\bf K})$
and each $y\in \bf K$, where $B(Y,y,r):=\{ z: z\in Y, d(y,z)\le r \} $
is the ball in an ultrametric space $Y$ with an ultrametric $d$ and a point
$y\in Y$. 
\par We have the following isomorphisms:
\par $L^g(T,{\cal R},\eta ,H)=L^g(T,{\cal R},\eta ,{\bf K})\otimes H$ and
\par $C^0_b(T,H)=C^0_b(T,{\bf K})\otimes H$, moreover,
$L^g(T,{\cal R},\eta ,{\bf K})$ is isomorphic with $c_0(\beta _L,{\bf K})$
and $C^0_b(T,{\bf K})$ is isomorphic with $c_0(\beta _C,{\bf K})$,
where $\beta _L$ and $\beta _C$ are ordinals, since $\Gamma _{\bf K}$
is discrete (see Chapter 5 \cite{roo}).
The locally $\bf K$-convex space $H^T$ is isomorphic with
$Y_1\otimes H$, where $Y_1:={\bf K}^T$. On the other hand,
the Banach space $c_0(\alpha ,{\bf K})$ in its weak topology
$\tau _w$ is isomorphic with ${\bf K}^T$ for $card (\alpha )=card (T)$,
since the valuation group of $\bf K$ is discrete in $(0,\infty )$
(see \S 8.203 in \cite{nari}).
Therefore, the ring of clopen subsets in $(c_0(\alpha ,{\bf K}),\tau _w)$
supplies $Y_1$ with the covering separating ring.
If $\mu _1$ and $\mu _2$ are $\bf K_s$-valued measures on Banach spaces
$Y_1$ and $Y_2$ with covering rings ${\cal R}_1$ and ${\cal R}_2$ respectively, 
then $\mu _1\otimes \mu _2$ is the $\bf K_s$-valued measure on $(Y_1\otimes Y_2, 
{\cal R}_1\times {\cal R}_2)$. In the Banach space $c_0(\beta ,{\bf K})$
there exists the canonical base $(e_j: j\in \beta )$, where
$e_j:=(0,...,0,1,0,...)$ with $1$ on the $j$-th place.
With this standard base are associated projections $\pi _{j_1,...,j_n}(x):=
\sum_{l=1}^nx^{j_l}e_{j_l}$ for each $j_1,...,j_n\in \beta $
and each $n\in \bf N$ and for each vector $x\in c_0(\beta ,{\bf K})$
with coordinates $x^j\in \bf K$ in the standard base.
Consider a covering separating ring $\cal R$ of $c_0(\beta ,{\bf K})$
such that 
\par ${\cal R}:=\bigcup_{j_1,...,j_n\in \beta ;n\in \bf N}
(\pi _{j_1,...,j_n})^{-1}(Bco(span_{\bf K}(e_{j_1},...,e_{j_n}))$,
\\
where $span_{\bf K}(z_l: l\in \gamma )$ $:=\{ x: x\in c_0(\beta ,{\bf K});$
$x=\sum_{j\in \zeta }a^jz_j;$ $a^j\in {\bf K};$
$card (\zeta )<\aleph _0 \} $ for each $\gamma \subset \beta $.
On the completion ${\cal R}_{\mu }$ there exists
a probability $\bf K_s$-valued measure $\mu $ generated by a 
bounded cylindrical distribution as in \S 2.8, \S 2.9 or \S 2.15.
For example, each $\mu _j(dx):=f_j(x)w(dx)$
is a measure on $\bf K$, where $f_j\in L({\bf K},
{\cal R}({\bf K}),w,{\bf K_s})$, $w$ is either the Haar measure or any
other probability measure on $\bf K$, $\mu _j=\pi _j(\mu )$
for each $j\in \beta $. 
\par In particular, for $card (\beta )\le \aleph _0$ and locally compact
$\bf K$ non-Archimedean infinite field with nontrivial valuation
there exists $\mu $ such that ${\cal R}_{\mu }\supset 
Bco(c_0(\beta ,{\bf K}))$.
For this consider on the Banach space $c_0:=c_0(\omega _0,{\bf K})$ 
a linear operator $J\in L_0(c_0)$, where $L_0(H)$ denotes
the Banach space of compact $\bf K$-linear operators on
the Banach space $H$, such that
$Je_i=v_ie_i$ with $v_i\ne 0$ for each $i$ and a measure $\nu (dx):=f(x)w(dx)$,
where $f: {\bf K}\to B({\bf K},0,r)$ with $r\ge 1$ 
is a function belonging to the space
$L({\bf K},{\cal R}_w,w,{\bf K_s})$ such that $\lim_{|x|\to \infty }f(x)=0$
and $\nu ({\bf K})=1$, $\| S\|_{\nu }>0$ for each clopen subset $S$ in $\bf K$,
for example, when $f(x)\ne 0$ $w$-almost everywhere. 
In particular we can choose $\nu $ with $\| \nu \|=1$.
In view of Lemma 2.3 and Theorem 2.30 from II \cite{lulapm} and
Theorem 2.8 there exists a product measure \\
$(i)\quad \mu (dx):=\prod_{i=1}^{\infty }\nu _i(dx^i)$ on the ring
$Bco(c_0)$ of clopen subsets of $c_0$, where 
\par $(ii)$ $\nu _i(dx^i):=f(x^i/v_i)\nu (dx^i/v_i)$.
\par Consider, for example, the particular case of $X=C^0(T,H)$
with compact $T$. If $t_0\in T$ is an isolated point, then
$C^0(T,H)=C^0(T\setminus \{ t_0 \},H)\oplus H$, so we consider the case
of $T$ dense in itself. Let $Z$ be a compact subset without
isolated points in a local field $\bf K$. Then 
the Banach space $C^0(Z,{\bf K})$ has the Amice polynomial 
orthonormal base
$Q_m(x)$, where $x\in Z$, $m\in {\bf N_o}:=\{ 0,1,2,... \} $ \cite{ami}.
Each $f\in C^0$ has a decomposition $f(x)=\sum_ma_m(f)Q_m(x)$
such that $\lim_{m\to \infty }a_m=0$, where $a_m\in \bf K$.
These decompositions establish the isometric isomorphism
$\theta : C^0(T,{\bf K})\to c_0(\omega _0,{\bf K})$
such that $\| f\|_{C^0}=\max_m|a_m(f)|=\| \theta (f)\|_{c_0}$.
If $u_i$ are roots of basic polynomils $Q_m$ as in \cite{ami},
then $Q_m(u_i)=0$ for each $m>i$. The set $\{ u_i: i \} $ is dense in $T$. 
\par The locally $\bf K$-convex space
$X=X(T,H)$ is isomorphic with the tensor product 
$X(T,{\bf K})\otimes H$ (see \S 4.R \cite{roo} and \cite{nari}).
If $J_i\in L_0(Y_i)$ is nondegenerate for each $i=1,2$, that is, 
$ker (J_i)=\{ 0 \} $, then $J:=J_1\otimes J_2\in L_0(Y_1\otimes Y_2)$
is nondegenerate (see also Theorem 4.33 \cite{roo}).
If $X(T,{\bf K})$ and $H$ are of separable type over 
a non-Archimedean locally compact infinite field $\bf K$ with 
nontrivial valuation, then 
we can construct a measure $\mu $ on $X$ such that ${\cal R}_{\mu }\supset
Bco(X)$. The case $H^T$ we reduce to $(c_0(\alpha ,{\bf K}), \tau _w)\otimes
H$ as above. Put 
$Y_1:=X(T,{\bf K})$ and $Y_2:=H$ 
and $J:=J_1\otimes J_2\in L_0(Y_1\otimes Y_2),$
where $J_1e_m:=\alpha _me_m$ such that $\alpha _m\ne 0$ for each $m$
and $\lim_i \alpha _i=0$. Take $J_2$ also nondegenerate. Then $J$
induces a product measure $\mu $ on $X(T,H)$
such that $\mu =\mu _1\otimes \mu _2$, where $\mu _i$ are
measures on $Y_i$ induced by $J_i$ due to Formulas $(i,ii)$.
Analogously considering the following subspace
$X_0(T,H)$and operators $J:=J_1\otimes J_2\in L_0(X_0(T,{\bf K})\otimes H)$
we get the measures $\mu $ on it also, where $t_0\in T$
is a marked point. 
On the other hand, the space $X(T,H)$ is Lindel\"of (see \S 3.8
\cite{eng}), hence each subset $U$ open in $X(T,H)$ is a countable
union of clopen subsets. Hence the characteristic function
$Ch_U$ of $U$ belongs to $L(X,{\cal R},\mu ,{\bf K_s})$, since
$\| \mu \| =\sup_xN_{\mu }(x)<\infty $, consequently, ${\cal R}_{\mu }
\ni U$.
In general the conditon ${\cal R}_{\mu }\supset Bco(X_0)$ is not imposed in
\S 4.1 and in Definiton 4.2, so ${\cal R}_{\mu }$ may be any covering ring
of $X_0$.
\par For each finite number of pairwise distinct 
points $(t_0,t_1,...,t_n)$ in $T$ and points
$(0,z_1,...,z_n)$ in $H$ there exists a closed subset
\\
$X(T,H;(t_0,t_1,...,t_n);(0,z_1,...,z_n))$ in $X(T,H)$ such that
\\
$X(T,H;(t_0,t_1,...,t_n);(0,z_1,...,z_n))=(0,z_1,...,z_n)+
X(T,H;(t_0,t_1,...,t_n);(0,...,0))$, \\ 
where
$X(T,H;(t_0,t_1,...,t_n);(0,...,0))$ is the $\bf K$-linear subspace 
in $X(T,H)$. Therefore, 
\par $(iii)$ rings $F_{t_2,t_1}^{-1}({\cal R}(H))$ and
$F_{t_4,t_3}^{-1}({\cal R}(H))$ are independent subrings
in the ring ${\cal R}(X(T,H))$, when $(t_1,t_2)$ and 
$(t_3,t_4)$ satisfy Condition $4.2.(i)$, where
covering rings ${\cal R}(H)$ and ${\cal R}(X)$ and measures
$\mu $, $\mu _1$ and $\mu _2$ are as above.
\par Put $P(t_1,x_1,t_2,A):=\mu ( \{ f: f\in X(T,H;(t_0,t_1);
(0,x_1)),\quad f(t_2)\in A \} )$
for each $t_1\ne t_2\in T,$ $x_1\in H$ and $A\in {\cal R}(H)$. In view of $(iii)$
we get, that $P$ satisfies Conditions $3.2.(i-iv).$ By the above construction
(and Proposition 3.3.1 also) the Markov cylindrical distribution
${\tilde \mu }_{x_0}$
induced by $\mu $ is bounded, since $\mu $ is bounded, where
$x_0=0$ for $X_0(T,H)$. Let $\Upsilon $ be a set of elementary events
$\omega :=\{ f: f\in X(T,H;(t_0,t_1,...,t_n);(0,x_1,...,x_n)) \} ,$
where $\Lambda _{\omega }$ is a finite subset of $\bf N$,
$x_i\in H$, $(t_i: i\in \Lambda _{\omega } )$ is a subset of $T\setminus 
\{ t_0\}$
of pairwise distinct points. There exists the ring $\tilde {\sf U}$
of cylindrical subsets of $X_0(T,H)$ induced by projections
$\pi _s: X_0(T,H)\to H^s,$ where $H^s:=\prod_{t\in s}H_t,$
$s=(t_1,...,t_n)$ are finite subsets of $T$, $H_t=H$ for each $t\in T$.
This induces the covering ring ${\cal R}({\Upsilon })$ of
${\Upsilon }$, where $({\Upsilon },{\cal R}({\Upsilon }),\nu )$
is the image of $(X_0(T,H),{\tilde {\sf U}},{\tilde \mu }_{x_0})$
due to Proposition 2.6.
In view of the Kolmogorov theorem 2.15.2 and \S 2.5 $\tilde \mu _{x_0}$
on $(X_0(T,H),{\tilde {\sf U}})$ induces the probability measure
$\nu $ on $(\Upsilon ,{\cal R}_{\nu }(\Upsilon ))$.
For each probability space $(\Omega ,{\sf F},\lambda )$
and a measurable mapping $\pi : \Omega \to \Upsilon $,
that is, $\pi ^{-1}({\cal R}(\Upsilon ))\subset {\sf F}$,
such that $\pi (\lambda )=\nu $ we get the space $L^q(\Omega ,{\sf F},\lambda ,
X_0)$ and the realization of the stochastic process $\xi (t,\omega )$.
In particular, we can take $\Omega =\Upsilon $ and $\pi =id$.
\par In the case $X(T,H)=H^T$ (apart from $C^0(T,H)$ and
$L^q(T,{\cal R},\eta ,H)$) it is sufficient to take
any bounded linear operator $J_1$ on $Y_1$, that is, $J_1\in 
L(Y_1)$, that brings the difference, when $card (T)\ge \aleph _0$.
\par Therefore, using cylindrical distributions we get examples of such 
measures $\mu $ for which stochastic processes exist.
Hence to each such measure
on $X_0(T,H)$ there corresponds the stochastic process.
\par Evidently, on the field $\bf K$ there exists a family $\Psi _{\bf K}$
of inequivalent $\bf K_s$-valued measures of the cardinality
$card (\Psi _{\bf K})\ge card ({\bf K})$, since the 
subfamily of atomic measures satisfies this inequality.
In the particular case of ${\bf C_p}\supset
{\bf F}\supset \bf Q_p$ or ${\bf F}=
\bf F_p(\theta )$ we use the Haar measure $w$ also for which
$card (L^q({\bf F},w,Bco({\bf F}),{\bf K_s}))={\sf c}:=card ({\bf R})$
for each $1\le q\le \infty $.
In view of the non-Archimedean variant of the Kakutani theorem
(see II.3.5 \cite{lulapm}) we get the inequalities
for $card (\Psi )$, since $card (H)=card (\beta _H)card ({\bf K})$.
In particular, for $T=B({\bf K},0,r)$ with $r>0$ and a locally compact field
$\bf K$ either ${\bf K}\supset \bf Q_p$ or ${\bf K}=\bf F_p(\theta )$
considering all operators
$J:=J_1\otimes J_2\in L_0(Y_1\otimes Y_2)$ and the corresponding 
measures as above we get ${\sf c}^{\aleph _0}=\sf c$ 
inequivalent measures for each chosen $f$.
\par {\bf Note.} Evidently, this theorem is also true for $C^0(T,H)$, that
follows from the proof. If take $\nu $ with $supp (\nu )=B({\bf K},0,1),$
then repeating the proof it is possible to construct $\mu $ with
$supp ( \mu )\subset B(C^0(T,{\bf K}),0,1)\times B(H,0,1).$
Certainly such measure $\mu $ 
can not be quasi-invariant relative
to shifts from a dense $\bf K$-linear subspace in $C^0(T,H)$,
but (starting from the Haar measure $w$ on $\bf F$)
$\mu $ can be constructed quasi-invariant relative to a dense
additive subgroup $G'$ of $B(C^0(T,{\bf K}),0,1)\times B(H,0,1)$, 
moreover, there exists $\mu $ for which $G'$ is also
$B({\bf K},0,1)$-absolutely convex.
\section{Poisson processes.}
{\bf 5.1. Definition.} Let $T$ be an additive group such that
$T\subset B({\bf K_s},0,r)$ and $0\ne \rho \in \bf K_s$
with $|\rho |r<s^{1/(1-s)}$, where ${\bf K_s}$ is a field such that
${\bf Q_s}\subset {\bf K_s}\subset \bf C_s$, ${\bf K_s}$
is complete as the ultrametric space.
Consider a stochastic process $\xi \in L^q(\Omega ,{\sf F},\lambda ,
X_0(T,H))$ such that the transition measure
has the form 
\par $P(t_1,x,t_2,A):=P(t_2-t_1,x,A):=Exp(-\rho (t_2-t_1))
P(A-x)$ \\
(see \S 3.2 and \S 4.2) for each $x\in H$ and $A\in {\cal R}_H$
and $t_1$ and $t_2$ in $T$,
where 
\par $Exp(x):=\sum_{n=0}^{\infty }x^n/n!$. 
Then such process is called the Poisson process.
\par {\bf 5.2. Proposition.} {\it Let
\par $P(A-x)=\int_HP(-x+dy)P(A-y)$ 
for each $x\in H$ and $A\in {\cal R}_H$ and 
\par $P(H)=1$
and $\| P\| =1$, 
\par then there exists
a measure $\mu $ on $X_0(T,H)$ for which the Poisson process
exists.}
\par {\bf Proof.} The exponential function
converges if $|x|<s^{1/(1-s)}$, since $|n!|_s^{-1}\le s^{(n-1)/(s-1)}$
for each $0<n\in \bf Z$ in accordance with Lemma 4.1.2 \cite{khipr}.
Hence $|Exp(-\rho (t_2-t_1))-1|<1$ for each $t_1$ and $t_2$ in $T$,
consequently, $|Exp(-\rho (t_2-t_1))|=1$.
Take 
\par $\mu _{t_1,...,t_n}:=P(t_2-t_1,0,*)...P(t_n-t_{n-1},0,*)$
on 
\par ${\cal R}_{t_1,...,t_n}:={\cal R}_{t_1}\times ...\times {\cal R}_{t_n}$
\\
for each pairwise distinct points $t_1,...,t_n\in T$,
where 
\par $\mu _{t_1,...,t_n}=\pi _{t_1,...,t_n}(\mu )$
and 
\par $\pi _{t_1,...,t_n}: X_0(T,H)\to H_{t_1}\times ...\times H_{t_n}$
\\
is the natural projection, $H_t=H$ for each $t\in T$.
There is a family $\Lambda $ of all finite subsets of $T$
directed by inclusion.
In view of Theorem 2.14 the cylindrical distribution
$\mu $ generated by the family $P(t_2-t_1,x,A)$ has an extension
to a measure on $X_0(T,H)$. All others conditions are satisfied 
in accordance with \S 4.2 and \S 5.1.
\par {\bf 5.3. Note.} 
Let $K$ be a complete ultrametric space with an ultrametric $d$,
that is,
\par  $d(x,y)\le \max (d(x,z), d(y,z))$ for each $x, y, z \in X$.
\\
Let 
\par $d(x,y):=\max_{1\le i \le n} d(x_i,y_i)$ \\
be the ultrametric in $K^n$,
where $x=(x_i: i=1,...,n)\in K^n$, $x_i\in K$.
Put 
\par ${\tilde K}^n:=(x\in K^n: x_i\ne x_j$
for each $i\ne j)$. \\
Supply ${\tilde K}^n$ with an ultrametric 
\par $\delta ^n_K(x,y):=d^n_K(x,y)/[\max (d^n_K(x,y),
d^n_K(x,({\tilde K}^n)^c), d(y,({\tilde K}^n)^c)]$,
\\
where $A^c:=K^n\setminus A$ for a subset $A\subset K^n$.
Then $({\tilde K}^n, \delta ^n_K)$ is the complete ultrametric
space. Let also $B^n_K$ denotes the collection of all
$n$-point subsets of $K$. Then the ultrametric 
$\delta ^n_K$ is equivalent with
the following ultrametric 
\par $d^{(n)}_K(\gamma ,\gamma '):=\inf_{\sigma \in
\Sigma _n}d^n_K((x_1,...,x_n), (x'_{\sigma (1)},...,x'_{\sigma (n)}))$,
\\
where $\Sigma _n$ is the symmetric group of $(1,...,n)$, $\sigma \in \Sigma _n$,
$\sigma : (1,...,n)\to (1,...,n)$; $\gamma , \gamma ' \in B^n_K$.
For each subset $A\subset K$ a number mapping $N_A: B^n_K\to \bf N_o$
is defined by the following formula: $N_A(\gamma ):=card(\gamma \cap A)$,
where ${\bf N}:=\{ 1,2,3,... \} $, ${\bf N_o}:=\{ 0,1,2,3,... \} $.
It remains to show, that $\delta ^n_K$ is the ultrametric for the ultrametric
space $(K,d)$.  For this we mention, that 
\par $(i)$ $\delta ^n_K(x,y)>0$, when
$x\ne y$, and $\delta ^n_K(x,x)=0$. 
\par $(ii)$ $\delta ^n_K(x,y)=\delta ^n_K(y,x)$,
since this symmetry is true for $d^n_K$ and for $[*]$ in the denumerator
in the formula defining $\delta ^n_K$.  To prove 
\par $(iii)$ $\delta ^n_K(x,y)\le
\max (\delta ^n_K(x,z), \delta ^n_K(z,y))$ \\
we consider the case
$\delta ^n_K(x,z)\ge \delta ^n_K(y,z)$, hence it is sufficient to show, that
$\delta ^n_K(x,y)\le \delta ^n_K(x,z)$.
Let 
\par $(a)$ $d^n_K(x,z)\ge \max (d^n_K(z,({\tilde K}^n)^c), d^n_K(x,({\tilde
K}^n)^c))$, \\
then $\delta ^n_K(x,z)=1$, hence $\delta ^n_K(x,y)\le \delta
^n_K(x,z)$, since  $\delta ^n_K(x,y)\le 1$
for each $x, y\in {\tilde K}^n$. Let 
\par $(b)$ $d^n_K(x,({\tilde K}^n)^c)>
\max (d^n_K(x,z), d^n_K(z,({\tilde K}^n)^c))$, then 
\par $\delta ^n_K(x,z)=
d^n_K(x,z)/d^n_K(x,({\tilde K}^n)^c)\le 1$.
Since $d^n_K(z,A):=\inf_{a\in A}d^n_K(z,a)$, then 
\par $d^n_K(z,({\tilde K}^n)^c)
\le \max (d^n_K(y,({\tilde K}^n)^c), d^n_K(y,z))$. \\
If $d^n_K(x,z)<
d^n_K(z,({\tilde K}^n)^c)$ and $d^n_K(x,y)\le d^n_K(x,z)$, then
\par $d^n_K(z,({\tilde K}^n)^c)\le d^n_K(x,({\tilde K}^n)^c)$. Hence
\par $d^n_K(x,y) \max (d^n_K(x,z), d^n_K(x, ({\tilde K}^n)^c), d^n_K(z,({\tilde
K}^n)^c))$ 
\par $\le d^n_K(x,z)\max (d^n_K(x,y), d^n_K(x,({\tilde K}^n)^c),
d^n_K(y,({\tilde K}^n)^c))$. \\
With the help of $(ii)$ the remaining cases
may be lightly written. 
\par {\bf 5.4. Notes and definitions.} As usually
let 
\par $B_K:=\bigoplus_{n=0}^{\infty }B^n_K$,
\\
where $B^0_K:=\{ \emptyset \} $ is a singleton,
$B_K\ni x=(x_n: x_n\in B^n_K, n=0,1,2,... )$. 
If a complete ultrametric space $X$ is not compact,
then there exists an increasing sequence of subsets $K_n\subset X$
such that $X=\bigcup_nK_n$ and $K_n$ are complete spaces in the induced
uniformity from $X$. Moreover, $K_n$ can be chosen clopen in $X$. 
Then the following space 
\par $\Gamma _X:=
\{ \gamma : \gamma \subset X$ and $card(\gamma \cap K_n)<\infty $
for each $n \} $ \\
is called the configuration space and it is isomorphic
with the projective limit $pr-\lim \{ B_{K_n}, \pi ^n_m, {\bf N} \} $,
where $\pi ^n_m(\gamma _m)=\gamma _n$ for each $m>n$ and $\gamma _n
\in B_{K_n}$. If $d_n$ denotes the ultrametric in $B_{K_n}$, then
$d_{n+1}|_{B_{K_n}}=d_n$, since $K_n\subset K_{n+1}$.
Then $\prod_{n=1}^{\infty }B_{K_n}=:Y$ in the Tychonoff product
topology is ultrametrizable, that induces the ultrametric in
$\Gamma _X$, for example, 
\par $\rho (x,y):=d_n(x_n,y_n)p^{-n}$ is the ultrametric in $\Gamma _X$,
\\
where $n=n(x,y):=\min_{(x_j\ne y_j)}j$,
$x=(x_j: j\in {\bf N}, x_j\in B_{K_j} ) $, $1<p\in \bf N$.
\par  Let $K\in \{ K_n: n\in {\bf N} \} $, then $m_K$
denotes the restriction $m|_K$, where $m: {\cal R}\to \bf K_s$
is a measure on a covering ring ${\cal R}_m$ of $X$, $K_n\in {\cal R}_m$
for each $n\in \bf N$. Suppose that $K^n_l\in {\cal R}_{m^n}$
for each $n$ and $l$ in $\bf N$, where ${\cal R}_m$ is the completion
of the covering ring ${\cal R}^n$ of $X^n$ relative to the product measure
$m^n=\bigotimes_{j=1}^nm_j$, $m_j=m$ for each $j$.
Then $m^n_K:=\bigotimes_{j=1}^n(m_K)_j$ is a measure
on $K^n$ and hence on ${\tilde K}^n$, when $m$ is 
such that $\| m|_{(K^n\setminus {\tilde K}^n)} \| =0$, for example,
non-atomic $m$, where $(m_K)_j=m_K$ for each $j$. 
Let $m(K_l)\ne 0$ for each $l\in \bf N$, $m(X)\ne 0$ and
$\| m\| <s^{1/(1-s)}$. Therefore, 
\par $(i)$ $P_{K,m}:=Exp(-m(K))\sum_{n=0}
^{\infty }m_{K,n}/n!$ \\
is a measure on ${\cal R}(B_K)$, where
\par ${\cal R}(B_K)=B_K\cap (\bigoplus_{n=0}^{\infty }{\cal R}_{m^n}),$
\\
$m_{K,0}$ is a probability measure on the singleton $B^0_K$, and
$m_{K,n}$ are images of $m^n_K$ under the following mappings:
\par $p^n_K: (x_1,...,x_n)\in {\tilde K}^n\to \{ x_1,...,x_n \} \in B^n_K$.
\\
Such system of measures
$P_{K,n}$ is consistent, that is, 
\par $\pi ^n_l(P_{K_l,m})=P_{K_n,m}$ for each
$n\le l$. \\
This defines the unique measure $P_m$ on 
${\cal R}(\Gamma _X)$, which is
called the Poisson measure, where
$\pi _n: Y\to B_{K_n}$ is the natural projection for each
$n\in \bf N$ (see for comparison the case of real-vauled
Poisson measures in \cite{shim}).
For each $n_1,..,n_l\in \bf N_o$ and disjoint subsets
$B_1,...,B_l$ in $X$ belonging to ${\cal R}_m$
there is the following equality:
\par $(ii)$ $P_m(\bigcap_{j=1}^l \{ \gamma : card(\gamma \cap B_i)=n_i \}
)=\prod_{i=1}^lm(B_i)^{n_i}Exp(-m(B_i))/n_i!$.
\par There exists the following embedding $\Gamma _X\hookrightarrow
S_X$, where 
\par $S_X:=\lim \{ E_{K_n}, \pi ^n_m, {\bf N} \} $ 
is the limit of an inverse mapping sequence,
\par $E_K:=\bigoplus_{l=0}
^{\infty }K^l$ for each $K\in \{ K_n: n=0,1,2,... \} $. \\
The Poisson measure $P_m$ on ${\cal R}(\Gamma _X)$ considered
above has an extension on ${\cal R}(S_X)$ such that
$\| P_m|_{S_X\setminus \Gamma _X}\| =0$. If each $K_n$ is a 
complete $\bf K$-linear space, then $E_K$ and $S_X$ are complete 
$\bf K$-linear spaces, since 
\par $S_X\subset (\prod_{n=1}^{\infty }
E_{K_n})$. \\
Then on ${\cal R}(S_X)$ there exists a Poisson measure
$P_m$, but without the restriction 
$\| m^n_K|_{K^n\setminus {\tilde K}^n}\| =0$, where
\par $(iii)$ $P_{K,m}:=Exp(-m(K))\sum_{n=0}^{\infty }
m^n_K/n!$, 
\par $\pi ^n_l(P_{K_l,m})=P_{K_n,m}$ for each $n\le l$.
\par {\bf 5.5. Corollary.} {\it Let suppositions of Proposition
5.2 be satisfied with $H=S_X$ for a complete $\bf K$-linear space
$X$ and $P(A)=P_m(A)$ for each $A\in {\cal R}(S_X)$, then there exists
a measure $\mu $ on $X_0(T,H)$ for which the Poisson process
exists.}
\par {\bf 5.6. Definition.} The stochastic process of Corollary 5.5
is called the Poisson process with values in $X$.
\par {\bf 5.7. Note.} If $\xi \in L^q(\Omega , {\sf F}, \lambda ;
X_0(T,H))$ is a stochastic process, then its mean value at the moment
$t\in T$ is defined by the following formula:
$$(i)\quad M_t(\xi ):=\int_{\Omega }\xi (t,\omega  )\lambda (d\omega ).$$
Let $H=\bf K$ be a field, where ${\bf Q_p}\subset {\bf K}\subset \bf C_p$,
let also $\lambda $ be with values in $\bf  K$. Suppose that
\par $\| \{ \omega : |\xi (t,\omega )|>R \} \|_{\lambda }=0$ 
for each $t\in T$ for some $R>0$. \\
Let $\rho \in B({\bf K},0,c)$
and let $T$ be a subgroup of $B({\bf K},0,r)$, 
where $R\max (c, r)<p^{1/(1-p)}$,
${\bf K}$ is the locally compact field.
Let 
\par ${\tilde P}^1: C^0(B({\bf K},0,c),{\bf K})\to 
C^1(B({\bf K},0,c),{\bf K})$ \\
be an antiderivation
operator (see also \S \S 54, 80 \cite{sch1}).
\par {\bf 5.8. Theorem.} (Non-Archimedean analog of the L\`evy
theorem.) {\it Let $\psi $ be a continuously differentiable function, 
from $T$ into ${\bf K}$ belonging to ${\tilde P}^1(
C^0(B({\bf K},0,c),{\bf K}))$ and $\psi (0)=0$.
Then there exists a stochastic process such that 
\par $M_t(Exp(-\rho \xi (t,\omega )))=Exp(-t\psi (\rho ))$ \\
for each $t$ in $T$ and each $\rho \in B({\bf K},0,c)$.}
\par {\bf Proof.} For the construction of $\xi $
consider solution of the following equation 
$$M_t[Exp(-\rho \xi (t,\omega ))]=Exp (-t\psi (\rho )).$$
Then $e(t)=e(t-s)e(s)$ 
for each $t$ and $s$ in $T$ and each $\rho \in
B({\bf K},0,c)$, where 
\par $e_{\rho }(t):=e(t):=M_t(Exp(-\rho \xi (t,\omega ))).$ \\
Hence 
\par $\partial e_{\rho }(t)/\partial \rho =-t\psi '(\rho )Exp(-t\psi
(\rho ))$, consequently, 
\par $\psi '(\rho )=t^{-1}\int_{\bf K}l
\quad EXP(-\rho l)P(\{ \omega : \xi (t,\omega )\in dl \} )$ \\
for each $t\ne 0$, where $EXP$ is the locally analytic extension of
$Exp$ on $\bf K$ (with values in $ \{ x: x\in {\bf C_p},\quad
|x-1|<1 \} $, see \cite{sch1}). In particular, 
\par $\lim_{t\to 0, t\ne 0}
t^{-1}\int_{\bf K}l\quad EXP(-\rho l)P(\{ \omega : \xi (t,\omega )
\in dl \} )$. \\
By the conditions of this theorem we have
\par  $\psi (\rho )=
{\tilde P}^1_{\beta }\psi '(\beta )|_0^{\rho }$. \\
Consider a measure $m$ on a separating covering ring
${\cal R}({\bf K})$ such that ${\cal R}({\bf K})\supset
Bco({\bf K})\cup \{ 0 \} $ with values in $\bf K$ such that
\par $m(dl):=\lim_{t\to 0, t\ne 0}lP(\{ \omega : \xi (t,\omega )
\in dl \} )/t$. \\
Therefore, 
\par $\psi (\rho )={\tilde P}^1_{\beta }(\int_{\bf K}EXP(-\beta l)m(dl))
|_0^{\rho }.$ \\
From $\psi (0)=0$ we have $e_{\rho }(1)=1$ for each $\rho $, consequently,
\par $\psi (\rho )=\rho m_0+\int_{\bf K}[1-EXP(-\rho l)]l^{-1}m(dl)$, \\
where $m_0:=m(\{ 0 \} )$, since
\par  $\lim_{l \to 0, l \ne 0}
[1-EXP(-\rho l)]/l=\rho $ and 
$$\lim_{\rho \to 0, \rho \ne 0}
\int_{B({\bf K},0,k)}[1-EXP(-\rho l)]l^{-1}m(dl)=0$$
for each $k>0$. Define a measure $n(dl)$ such that $n(\{ 0\} )=0$
and $n(dl)=l^{-1}m(dl)$ on ${\bf K}\setminus \{ 0 \} $, then
\par $\psi (\rho )=\rho m_0+\int_{\bf K}[1-EXP(-\rho l)]n(dl)$.
We search a solution of the problem in the form
\par $\xi (t,\omega )=tm_0+\int_{\bf K}l{\cal \eta }(t,dl,\omega )$, \\
where ${\cal \eta }(t,dl,\omega )$ is the measure on ${\cal R}({\bf K})$
for each fixed $t\in T$ and $\omega \in \Omega $ such that
its moments satisfy the Poisson distribution with the Poisson measure
$P_{tn}$, that is, 
\par $M_t[{\cal \eta }^k(t,dl,\omega )]=
\sum_{s\le k}a_{s,k}(tn)^s(dl)/s!$ \\
for each $t\in T$, where
$a_{0,j}=0$, $a_{1,j}=1$ and recurrently
$$a_{k,j}=k^j-\sum_{s=1}^k{k\choose s}a_{k-s,j}$$
for each $k\le j$, in particular, $a_{j,j}=j!$, that is,
$$a_{k,j}=\sum_{s_1+...+s_k=j, s_1\ge 1,...,s_k\ge 1}[j!/(s_1!...s_k!)].$$
Using the fact that the set of step functions is dense
in $L({\bf K},{\cal R}({\bf K}),n,{\bf C_p})$
we get
$$M_t[EXP(-\rho \int_{\bf K}l{\cal \eta }(t,dl,\omega ))]
=\lim_{\cal Z}M_t[\prod_jEXP(-\rho l_j{\cal \eta }(t,\delta _j,\omega ))]$$
$$=\lim_{\cal Z}\prod_jM_t[EXP(-\rho l_j{\cal \eta }(t,\delta _j,\omega ))]
=\lim_{\cal Z}EXP(-\rho t\sum_j(1-EXP(-\rho l_j))n(\delta _j))$$
$$=EXP[-\rho t\int_{\bf K}(1-EXP(-\rho l)n(dl)],$$
where ${\cal Z}$ is an ordered
family of partitions ${\cal U}$ of $\bf K$ into disjoint union of elements
of ${\cal R}({\bf  K})$, ${\cal U}\le \cal V$ in ${\cal Z}$
if and only if each element of the disjoint covering ${\cal U}$
is a union of elements of $\cal V$, $l_j\in \delta _j\in {\cal U}\in \cal Z$.
The limit
\par $\lim_{{\cal U}\in \cal Z}f({\cal U})=:\lim_{\cal Z}f=a$ \\
means that for each $\epsilon >0$ there exists ${\cal U}$
such that for each $\cal V$ with ${\cal U}\le \cal V$
we have 
\par $|a-f({\cal V})|<\epsilon $, \\
where $f({\cal U})$
is one of the functions defined as above with $l_j\in \delta _j\in {\cal U}$,
that is, 
\par $f({\cal U})=M_t[g\circ h(\eta )]$, \\
where $g\circ h(\eta )$ is the composition
of the continuous function $g$ and of 
\par $h(\eta )=\int_{\bf K}\zeta (y){\cal \eta }(t,dy,\omega )$ \\
with the step function $\zeta $.
We get the equation 
$$M_t[EXP(-\rho \xi (t,\omega )]=EXP(-\rho tm_0)
M_t[EXP(-\rho \int_{\bf K}l{\cal \eta }(t,dl,\omega ))].$$
In view of Corollary 5.5 it
defines the stochastic process with the probability space
$(\Omega ,{\sf F},\lambda )$, the existence of which follows
from the second half of \S 4.3.
\par {\bf 5.9. Note.} From the preceding results it follows, that
there are several specific features of non-Archimedean
stochastic processes and in particular Poisson processes
in comparison with the classical case. For this there are several reasons.
The non-Archimedean infinite field $\bf K$ with nontrivial valuation
has not any linear ordering compatible with its field structure.
In the non-Archimedean case there is not any indefinite integral.
Theory of analytic functions and elements has many specific features
in the non-Archimedean case \cite{esc,sch1}.
Moreover, interpretations of probabilities also are different 
\cite{khrif,khipr}. \\
\par We have started our collaboration with the investigation
of one problem formulated few years ago by A. Khrennikov and 
A.C.M. van Rooij. It was in the study of non-Archimedean analogs
of the Kolmogorov theorem for measures with values in non-Archimedean fields.
\par S. Ludkovsky is sincerely grateful to A. Khrennikov for his hospitality
at International Center for Mathematical Modeling of V\"axj\"o University.

\par Addresses: Theoretical Department, Institute of General Physics,
\par Russian Academy of Sciences,
\par Str. Vavilov 38, Moscow, 119991 GSP-1, Russia and
\par School of Mathematics and System Engineering,
\par International Center for Mathematical Modeling,
\par V\"axj\"o Universitet,
\par SE-351 95, Vejdes Plats 7, Sweden.
\end{document}